\documentclass[a4paper]{amsart}
\usepackage{amssymb} 
\usepackage[T1]{fontenc}
\usepackage[latin1]{inputenc}
\usepackage{amsfonts}
\usepackage{amsxtra}
\usepackage{ae}
\usepackage{slashed}
\usepackage{pdfsync}
\usepackage[all]{xy}
\usepackage{enumerate}
\usepackage{verbatim}
\usepackage{pgf,tikz}
\usetikzlibrary{arrows,patterns}
\include{diagram}

\newcommand*{\ket}{\rangle}
\newcommand*{\bra}{\langle}
\newcommand*{\ad}{\mathsf{ad}}

\newcommand*{\M}{\mathcal{M}}

\newcommand*{\E}{\mathcal{E}}
\newcommand*{\F}{\mathcal{F}}
\newcommand*{\G}{\mathcal{G}}
\renewcommand*{\H}{\mathcal{H}}

\newcommand*{\K}{\mathcal{K}}
\renewcommand*{\L}{\mathcal{L}}

\renewcommand*{\S}{\mathcal{S}}

\newcommand*{\V}{\mathcal{V}}

\renewcommand*{\max}{\mathsf{f}}
\newcommand*{\red}{\mathsf{r}}
\newcommand*{\Irr}{\mathsf{Irr}}

\newcommand*{\hit}{\rightharpoonup}
\newcommand*{\hitby}{\leftharpoonup}

\newcommand*{\KH}{\mathbb{K}}
\newcommand*{\LH}{\mathbb{L}}
\newcommand*{\CF}{\mathfrak{C}}
\newcommand*{\DF}{\mathfrak{D}}

\newcommand{\roots}{\mathbf{Q}}

\newcommand{\weights}{\mathbf{P}}
 
\renewcommand*{\top}{\mathsf{top}}

\DeclareMathOperator{\End}{End}
\DeclareMathOperator{\tr}{tr}

\DeclareMathOperator{\Hom}{Hom}

\DeclareMathOperator{\ind}{ind}

\DeclareMathOperator{\id}{id}

\newenvironment{bnum}
{\begin{list}{}
    {\setlength{\labelwidth}{15pt}
     \setlength{\leftmargin}{\labelwidth}
    }
}
{\end{list}}

\numberwithin{equation}{section}
\theoremstyle{change}
\newtheorem{theorem}{Theorem}[section]
\newtheorem{prop}[theorem]{Proposition}
\newtheorem{lemma}[theorem]{Lemma}

\newtheorem{definition}[theorem]{Definition}

\begin{document}

\title[Complex quantum groups and the assembly map]{Complex quantum groups and a deformation of the Baum-Connes assembly map}

\author{Andrew Monk}
\address{School of Mathematics and Statistics \\
         University of Glasgow \\
         University Place \\
         Glasgow G12 8SQ \\
         United Kingdom 
}

\email{a.monk.1@research.gla.ac.uk}

\author{Christian Voigt}
\address{School of Mathematics and Statistics \\
         University of Glasgow \\
         University Place \\
         Glasgow G12 8SQ \\
         United Kingdom 
}

\email{christian.voigt@glasgow.ac.uk}


\thanks{This work was supported by EPSRC grant no EP/K032208/1.
The second author was supported by the Polish National Science Centre grant no. 2012/06/M/ST1/00169.}


\maketitle

\begin{abstract}
We define and study an analogue of the Baum-Connes assembly map for complex semisimple quantum groups, that is, Drinfeld doubles 
of $ q $-deformations of compact semisimple Lie groups. 

Our starting point is the deformation picture of the Baum-Connes assembly map for a complex semisimple Lie group $ G $, which 
allows one to express the $ K $-theory of the reduced group $ C^* $-algebra of $ G $ in terms of the $ K $-theory of its associated Cartan motion group. 
The latter can be identified with the semidirect product of the maximal compact subgroup $ K $ acting on $ \mathfrak{k}^* $ via the coadjoint action. 

In the quantum case the role of the Cartan motion group is played by the Drinfeld double of the classical group $ K $, 
whose associated group $ C^* $-algebra is the crossed product of $ C(K) $ with respect to the adjoint action of $ K $.
Our quantum assembly map is obtained by varying the deformation parameter in the Drinfeld double construction applied to the standard deformation $ K_q $ 
of $ K $. We prove that the quantum assembly map is an isomorphism, thus providing a description of the $ K $-theory of complex quantum groups in terms of 
classical topology. 

Moreover, we show that there is a continuous field of $ C^* $-algebras which encodes both the quantum and classical assembly maps as well as a natural deformation 
between them. It follows in particular that the quantum assembly map contains the classical Baum-Connes assembly map as a direct summand. 
\end{abstract}

\section{Introduction}

Let $ G $ be a second countable locally compact group. The Baum-Connes conjecture asserts that the assembly map 
$$ 
\mu: K^{\top}_*(G) \rightarrow K_*(C^*_\red(G))
$$
is an isomorphism \cite{BC}, \cite{BCH}. Here $ C^*_\red(G) $ denotes the reduced group $ C^* $-algebra of $ G $. 
The conjecture has been established for large classes of groups, in particular, it holds for Lie 
groups \cite{Kasparov2}, \cite{WassermannCK}, \cite{PPdiracprincipalseries}, \cite{Lafforguebc}, 
and in fact for arbitrary almost connected groups \cite{CENconneskasparov}. 
This part of the Baum-Connes conjecture is also known as the Connes-Kasparov conjecture. 

In the case that $ G $ is a Lie group with finitely many components there exists an alternative description of the Baum-Connes assembly 
map \cite{Connesbook}, \cite{BCH}. Namely, if $ K \subset G $ is a maximal compact subgroup and $ \mathfrak{k} \subset \mathfrak{g} $ the 
corresponding inclusion of Lie algebras, there exists a continuous bundle of groups over $ [0,1] $ 
deforming $ G $ into its associated Cartan motion group $ G_0 = K \ltimes \mathfrak{g}/\mathfrak{k} $. 
Moreover, this bundle is trivial outside $ 0 $, and the resulting induced map 
$$ 
K_*(C^*(G_0)) \rightarrow K_*(C^*_\red(G)) 
$$ 
in $ K $-theory identifies naturally with the Baum-Connes assembly map. We refer to \cite{Malowthesis} for a detailed discussion. 
The deformation picture of the Baum-Connes map can be viewed as a $ K $-theoretic version of the Mackey analogy \cite{Mackeyanalogy}. 

In the special case of complex semisimple groups, Higson obtained a proof of the Connes-Kasparov conjecture in the deformation picture
by analysing the structure of the $ C^* $-algebras $ C^*(G_0) $ and $ C^*_\red(G) $, using the representation theory of 
complex semisimple groups and the Mackey machine \cite{Higsonmackeyanalogyktheory}. 

In the present paper we shall study the $ K $-theory of complex semisimple quantum groups from a similar perspective. 
These quantum groups are obtained by applying the Drinfeld double construction to $ q $-deformations of compact semisimple Lie groups \cite{PWlorentz}. 
We shall determine the $ K $-theory of the associated reduced group $ C^* $-algebras in terms of the $ K $-theory of a quantum analogue of the 
Cartan motion group. As in \cite{Higsonmackeyanalogyktheory} we start from an explicit description of the reduced group $ C^* $-algebras. 
 
Our definition of the quantum assembly map is very natural from the point of view of deformation quantization. To explain this, recall that 
the maximal compact subgroup $ K $ of a complex semisimple group $ G $ admits an essentially unique standard deformation. 
The classical double of the corresponding Poisson-Lie group \cite{CPbook} naturally identifies with the group $ G $. 
Moreover, if we instead consider the trivial Poisson bracket, then the corresponding classical double 
is precisely the Cartan motion group $ G_0 $. The Baum-Connes deformation may thus be viewed as rescaling the standard 
Poisson structure on $ K $ to zero and considering the corresponding classical doubles. 
In the quantum case, we obtain our deformation in the same way by rescaling the deformation parameter $ q $ to $ 1 $ and applying 
the Drinfeld double construction. Moreover we show that the various deformations appearing in this picture are all compatible in a natural sense. 

As already indicated above, our approach in the present paper relies on the representation theory of complex quantum 
groups \cite{VYcqg}, \cite{VYcqgplancherel}. It would be nice to find alternative proofs without having to invoke these results. Such an independent 
approach might shed some new light on the Baum-Connes isomorphism for classical complex semisimple groups as well. 

Let us now explain how the paper is organised. In section \ref{secprelim} we collect some preliminaries and 
fix our notation. Section \ref{secbc} contains a review of the deformation picture of the Baum-Connes assembly map for almost connected 
Lie groups, and more specifically complex semisimple groups. The left hand side of the assembly map is identified with 
the $ K $-theory of the Cartan motion group $ K \ltimes \mathfrak{k}^* $ associated to the complex group $ G $; here $ K $ 
is the maximal compact subgroup of $ K $ acting on the dual $ \mathfrak{k}^* $ of its Lie algebra $ \mathfrak{k} $ by the coadjoint
action. The main aim of section \ref{seccmg} is to explain how the Cartan motion group 
can be quantized in such a way that one obtains a continuous field of $ C^* $-algebras. 
Section \ref{seccg} deals with a corresponding quantization on the level of the complex group itself. Again we obtain 
a continuous field of $ C^* $-algebras. 
In section \ref{secqbc} we prove our first main result; more precisely, we define a quantum analogue of the Baum-Connes assembly map 
from the $ K $-theory of the quantum Cartan motion group with values in the $ K $-theory of the reduced group $ C^* $-algebra of 
the corresponding complex quantum group and show that this map is an isomorphism. 
Finally, in section \ref{secdefsquare} we bring together all the deformations studied in the previous sections. 
We introduce the concept of a deformation square, by which we mean a specific type of continuous field over the unit square. 
For each complex semisimple Lie group we construct a deformation square, 
which encodes both the quantum and classical assembly maps as well as a deformation between them. 

Let us conclude with some remarks on notation. The algebra of adjointable operators on a Hilbert space or Hilbert module $ \E $ is denoted 
by $ \LH(\E) $, and we write $ \KH(\E) $ for the algebra of compact operators. Depending on the context, the symbol $ \otimes $ denotes either 
the algebraic tensor product over the complex numbers, the tensor product of Hilbert spaces, or the minimal tensor product of $ C^* $-algebras. 

The authors would like to thank the Isaac Newton Institute for Mathematical Sciences, Cambridge, for support
and hospitality during the programme Operator Algebras: Subfactors and their applications,
where work on this paper was undertaken.

\section{Preliminaries} \label{secprelim}

In this section we review some background material on continuous fields of $ C^* $-algebras 
and on quantum groups.

\subsection{Continuous fields of $ C^* $-algebras and $ K $-theory} 

In this subsection we recall some basic definitions and facts on continuous fields of $ C^* $-algebras and their $ K $-theory, 
see \cite{Blancharddef}, \cite{Nilsenbundles}, \cite{Williamscrossedproducts}. 

If $ A $ is a $ C^* $-algebra we write $ ZM(A) $ for its central multiplier algebra, that is, the center of the multiplier 
algebra of $ A $. 

\begin{definition} 
Let $ X $ be a locally compact space. A $ C_0(X) $-algebra is a $ C^* $-algebra $ A $ together with a nondegenerate 
$ * $-homomorphism $ C_0(X) \rightarrow ZM(A) $. 
\end{definition} 

We will usually omit the homomorphism $ C_0(X) \rightarrow ZM(A) $ from our notation and simply write $ f a $ for the 
action of $ f \in C_0(X) $ on $ a \in A $. A morphism of $ C_0(X) $-algebras is a $ * $-homomorphism $ \phi: A \rightarrow B $ such 
that $ \phi(f a) = f \phi(a) $ for all $ a \in A $ and $ f \in C_0(X) $. In particular, two $ C_0(X) $-algebras are isomorphic
if there exists a $ * $-isomorphism between them which is compatible with the $ C_0(X) $-module structures. 

Assume $ A $ is a $ C_0(X) $-algebra and let $ x \in X $. The fibre of $ A $ at $ x $ is the $ C^* $-algebra 
$$
A_x = A/I_x A 
$$
where $ I_x \subset C_0(X) $ is the ideal of all functions vanishing at $ x $. We note that $ I_x A $ is a closed two-sided 
ideal of $ A $ due to Cohen's factorisation theorem. 
If $ a \in A $ we write $ a_x \in A_x $ for the image under the canonical projection homomorphism. 

\begin{definition} 
Let $ X $ be a locally compact space. A continuous field of $ C^* $-algebras over $ X $ is a $ C_0(X) $-algebra $ A $ 
such that for each $ a \in A $ the map $ x \mapsto \|a_x\| $ is continuous.  
\end{definition} 

A basic example of a continuous field over $ X $ is $ A = C_0(X, D) $ where $ D $ is an arbitrary $ C^* $-algebra and the 
action of $ C_0(X) $ is by pointwise multiplication. Continuous fields which are isomorphic to fields of this form are called \emph{trivial}. 

Let $ A $ be a $ C_0(X) $-algebra and assume that $ U \subset X $ is a closed subset. In the same way as in the 
construction of the fibre algebra $ A_x $ one constructs a $ C^* $-algebra $ A_U $ by setting 
$$
A_U = A/I_U A
$$
where $ I_U \subset C_0(X) $ is the ideal of all functions vanishing on $ U $. Note that we have $ (A_U)_V \cong A_V $ if $ V \subset U \subset X $ 
are closed subsets. In particular, $ A_U $ is a $ C_0(U) $-algebra with the same fibres as $ A $ at the points of $ U $. 
If $ A $ is a continuous field over $ X $ then $ A_U $ is a continuous field over $ U $. We call $ A_U $ the restriction of $ A $ to $ U $. 

Similarly, let $ A $ be a $ C_0(X) $-algebra and assume that $ U \subset X $ is an open subset. We let $ A_U $ be the kernel of the canonical 
projection $ A \rightarrow A_{X \setminus U} $. Equivalently, we may describe $ A_U $ as 
$$
A_U = C_0(U) A 
$$
in this case. 
In the same way as for closed subsets, we call $ A_U $ the restriction of $ A $ to $ U $. If $ A $ is a continuous field over $ X $ 
then $ A_U $ is a continuous field over $ U $, with the same fibres as $ A $ at the points of $ U $. 

A continuous field of $ C^* $-algebras $ A $ over a locally compact space $ X $ is called trivial away from $ x \in X $ if 
the restricted field $ A_{X \setminus \{x\}} $ is trivial. 

Let us now review a basic fact from the $ K $-theory of continuous fields which will play a central role in our considerations. 
Assume that $ A $ is a continuous field of $ C^* $-algebras over the unit interval $ X = [0, 1] $. If $ A $ is trivial away from $ x = 0 $, 
then we obtain a canonical induced homomorphism between the $ K $-groups of the fibres $ A_0 $ and $ A_1 $. Indeed, restriction to $ 0 $ yields a short 
exact sequence 
$$
\xymatrix{
0 \ar@{->}[r] & A_{(0,1]} \ar@{->}[r] & A \ar@{->}[r] & A_0 \ar@{->}[r] & 0
}
$$
of $ C^* $-algebras, and triviality of the field away from $ 0 $ means that the kernel in this sequence is isomorphic to $ C_0((0,1], A_1) $. 
In particular, the $ K $-theory of $ A_{(0,1]} $ vanishes, and therefore the $ 6 $-term exact sequence yields an isomorphism $ K_*(A) \rightarrow K_*(A_0) $. 
Combining this with the projection to the fibre at $ 1 $, the resulting diagram 
$$
\xymatrix{
K_*(A_0) \ar@{<-}[r]^{\;\; \cong} & K_*(A) \ar@{->}[r] & K_*(A_1) 
}
$$
yields a homomorphism $ i: K_*(A_0) \rightarrow K_*(A_1) $. This construction will be used repeatedly in the sequel. 

If the continuous field $ A $ is unital it is easy to describe the map $ i $ directly on the level of projections and unitaries. Since this will 
become important later on, let us briefly review this. 
Assume that $ p_0 \in M_n(A_0) $ is a projection. We can lift $ p_0 $ to a positive element $ q $ in $ M_n(A) $. Since $ p_0 - p_0^2 = 0 $ 
we see that $ \|q(\tau) - q(\tau)^2\| < \frac{1}{4} $ for $ \tau $ sufficiently small. In particular, since $ \frac{1}{2} $ is not in the spectrum of $ q(\tau) $ 
for such $ \tau $, functional calculus allows us to find a lift $ p \in M_n(A) $ of $ p_0 $ such that all elements $ p(\tau) \in M_n(A_\tau) $ are 
projections for $ \tau $ small. Since the field is trivial outside $ 0 $ we can in fact arrange that $ p(\tau) $ is a projection for all $ \tau \in [0,1] $. 
The image of $ [p_0] \in K_0(A_0) $ under the map $ i $ is then the class $ [p(1)] \in K_0(A_1) $. In particular, on the level of $ K $-theory the 
previous construction is independent of the choice of the lift $ p $. 

In a similar way, assume that $ u_0 \in M_n(A_0) $ is a unitary. Then we can lift $ u_0 $ to an element $ u $ in $ M_n(A) $, 
and for $ \tau $ small the resulting elements $ u(\tau) $ will satisfy $ \|u(\tau) u(\tau)^* - 1 \| < 1 $ and $ \|u(\tau)^* u(\tau) - 1 \| < 1 $. 
It follows that $ u(\tau) $ is invertible for $ t $ small, 
and by triviality of the field outside $ 0 $ we can assume that $ u(\tau) $ is 
invertible for all $ \tau \in [0, 1] $. Using polar decomposition we may as well arrange for $ u(\tau) $ to be unitary. In any case, 
the image of $ [u_0] \in K_1(A_0) $ under $ i $ is the class $ [u(1)] \in K_1(A_1) $. Again, upon passing to $ K $-theory this does not 
depend on the choice of the lifting section $ u $.

\subsection{Quantum groups} 

In this subsection we give a quick introduction to the theory of quantum groups with a particular focus on complex semisimple quantum groups. 
We refer to \cite{CPbook}, \cite{KS}, \cite{VYcqg} for more details. 

Let $ G $ be a simply connected complex semisimple Lie group and let $ \mathfrak{g} $ be its Lie algebra. 
We will fix a Cartan subalgebra $ \mathfrak{h} \subset \mathfrak{g} $ and a set $ \alpha_1, \dots, \alpha_N $ of simple roots. 
Moreover we let $ \mathfrak{k} \subset \mathfrak{g} $ be the Lie algebra of the maximal compact subgroup $ K \subset G $ with maximal torus $ T $ 
such that $ \mathfrak{t} \subset \mathfrak{h} $, where $ \mathfrak{t} $ is the Lie algebra of $ T $.  
Let $ ( \;, \;) $ be the bilinear form on $ \mathfrak{h}^* $ obtained by rescaling the Killing form such that $ (\alpha, \alpha) = 2 $ 
for the shortest root of $ \mathfrak{g} $, and set $ d_i = (\alpha_i, \alpha_i)/2 $. 
Moreover we let $ \roots \subset \weights \subset \mathfrak{h}^* $ be the root and weight lattices of $ \mathfrak{g} $, respectively. 
The set $ \weights^+ \subset \weights $ of dominant integral weights is the set of all non-negative integer combinations of the fundamental 
weights $ \varpi_1, \dots, \varpi_N $, the latter being defined by stipulating $ (\varpi_i, \alpha_j) = \delta_{ij} d_j $. Equivalently, 
the fundamental weights are the dual basis to the basis formed by the coroots $ \alpha_1^\vee, \dots, \alpha_N^\vee $, where $ \alpha_j^\vee = d_j^{-1} \alpha_j $. 

The $ C^* $-algebra $ C(K) $ of continuous functions on the maximal compact subgroup $ K $ can be obtained 
as a completion of the algebra of matrix coefficients of all finite dimensional representations of $ \mathfrak{g} $. 
In the same way one constructs a $ C^* $-algebra $ C(K_q) $ for $ q \in (0,1) $ as the completion of the algebra 
of matrix coefficients of all finite dimensional integrable representations of the quantized universal enveloping algebra $ U_q(\mathfrak{g}) $ 
associated with $ \mathfrak{g} $. Let us briefly review these constructions. 

We recall that $ U_q(\mathfrak{g}) $ is generated by elements $ K^q_\mu $ for $ \mu \in \weights $ and $ E^q_i, F^q_j $ for $ 1 \leq i,j \leq N $ 
satisfying a deformed version of the Serre presentation for $ U(\mathfrak{g}) $. 
More precisely, we will work with the Hopf algebra $ U_q(\mathfrak{g}) $ as defined in \cite{VYcqg}. We let $ U_1(\mathfrak{g}) = U(\mathfrak{g}) $ 
be the classical universal enveloping algebra of $ \mathfrak{g} $ and write
$ H^1_i, E^1_i, F^1_i $ for the Serre generators of $ U(\mathfrak{g}) $. 

In analogy with the classical case one has the notion of a weight module for $ U_q(\mathfrak{g}) $. 
Every finite dimensional weight module is completely reducible, and the irreducible finite dimensional weight modules of $ U_q(\mathfrak{g}) $ are 
parameterized by their highest weights in $ \weights^+ $ as in the classical theory. We will write $ V(\mu)^q $ for the module associated to $ \mu \in \weights^+ $ 
and $ \pi^q_\mu: U_q(\mathfrak{g}) \rightarrow \End(V(\mu)^q) $ for the corresponding representation. The direct sum of the maps $ \pi^q_\mu $ induces an embedding 
$ \pi^q: U_q(\mathfrak{g}) \rightarrow \prod_{\mu \in \weights^+} \End(V(\mu)^q) $. For $ q \in (0,1) $ we may use this to define $ H^q_j $ to be the unique 
element in $ \prod_{\mu \in \weights^+} \End(V(\mu)^q) $ such that $ K_j^q = q^{d_j H^q_j} $, where $ K_j^q = K^q_{\alpha_j} $. 

If we fix a $ * $-structure on $ U_q(\mathfrak{g}) $ as in \cite{VYcqg}, then the representations $ V(\mu)^q $ for $ \mu \in \weights^+ $ are unitarizable. 
In fact, one can identify the underlying Hilbert spaces $ V(\mu)^q $ with the Hilbert space $ V(\mu)^1 = V(\mu) $ of the corresponding representation of 
$ U(\mathfrak{g}) $, in such a way that the operators $ H^q_i, E^q_i, F^q_i $ define continuous families in $ \KH(V(\mu)) $ for all $ \mu \in \weights^+ $, 
compare \cite{NTKhomologydirac}. 

For $ q \in (0,1) $ we define $ \CF^\infty(K_q) $ to be the space of all matrix coefficients of finite dimensional weight modules over $ U_q(\mathfrak{g}) $. 
This space becomes a Hopf $ * $-algebra with multiplication, comultiplication, counit and antipode defined in such a way that the 
canonical evaluation pairing $ U_q(\mathfrak{g}) \times \CF^\infty(K_q) \rightarrow \mathbb{C} $ satisfies 
\begin{align*}
(XY, f) &= (X, f_{(1)}) (Y, f_{(2)}), \qquad (X, fg) = (X_{(2)}, f) (X_{(1)}, g) 
\end{align*}
and 
\begin{align*}
(\hat{S}(X), f) &= (X, S^{-1}(f)), \qquad (\hat{S}^{-1}(X), f) = (X, S(f)) 
\end{align*} 
\begin{align*}
(X, f^*) &= \overline{(\hat{S}^{-1}(X)^*, f)} 
\end{align*} 
for $ X, Y \in U_q(\mathfrak{g}) $ and $ f, g \in \CF^\infty(K_q) $. Here we use the Sweedler notation $ \hat{\Delta}(X) = X_{(1)} \otimes X_{(2)} $ 
for the coproduct of $ X \in U_q(\mathfrak{g}) $ and write $ \Delta(f) = f_{(1)} \otimes f_{(2)} $ for the coproduct of $ f \in \CF^\infty(K_q) $. Similarly,  
we denote by $ \hat{S}, \hat{\epsilon} $ the antipode and counit of $ U_q(\mathfrak{g}) $, and write $ S, \epsilon $ for the 
corresponding maps for $ \CF^\infty(K_q) $, respectively. 

By definition, the Hopf $ * $-algebra $ \CF^\infty(K_q) $ is the algebra of polynomial functions on the compact quantum group $ K_q $. 
We note that $ \CF^\infty(K_q) $ has a linear basis of matrix coefficients $ u^\mu_{ij} = \bra e^\mu_i|\bullet| e^\mu_j \ket $, 
where $ \mu \in \weights^+ = \Irr(K_q) $ is the set of equivalence classes of irreducible representations of $ K_q $
and $ e^\mu_1, \dots, e^\mu_n $ is an orthonormal basis of $ V(\mu) $. We will always consider matrix coefficients with respect to a basis of weight vectors. 
In terms of matrix elements, the normalized Haar functional $ \phi: \CF^\infty(K_q) \rightarrow \mathbb{C} $ is given by 
$$
\phi(u^\mu_{ij}) = 
\begin{cases}
1 & \text{if } \mu = 0 \\
0 & \text{else. }
\end{cases}
$$

Since $ K_q $ is a compact quantum group the algebra $ \CF^\infty(K_q) $ is an algebraic quantum group in the sense of van Daele \cite{vDadvances}. 
We shall write $ \DF(K_q) $ for its dual in the sense of algebraic quantum groups. Explicitly, the dual is is given by the algebraic direct sum 
$$
\DF(K_q) = \text{alg-}\bigoplus_{\mu \in \weights^+} \KH(V(\mu))  
$$
with the $ * $-structure arising from the $ C^* $-algebras $ \KH(V(\mu)) $. We denote by $ p_\eta $ the central projection in $ \DF(K_q) $ corresponding 
to the matrix block $ \KH(V(\eta)) $ for $ \eta \in \weights^+ $. 

The left and right Haar functionals for $ \DF(K_q) $ are given by 
$$
\hat{\phi}(x) = \sum_{\mu \in \weights^+} \dim_q(V(\mu)) \tr(K_{2\rho} p_\mu x), 
\qquad \hat{\psi}(x) = \sum_{\mu \in \weights^+} \dim_q(V(\mu)) \tr(K_{-2 \rho} p_\mu x), 
$$
respectively. Here $ \dim_q(V(\mu)) $ denotes the quantum dimension of $ V(\mu) $, and $ K_{2\rho} \in U_q(\mathfrak{g}) $ is the generator 
associated with $ 2 \rho $, where $ \rho $ is the half-sum of all positive roots.  

Given the basis of matrix coefficients $ u^\mu_{ij} $ in $ \CF^\infty(K_q) $ we obtain a linear basis of matrix units $ \omega^\mu_{ij} $ of $ \DF(K_q) $ satisfying 
$$
(\omega^\mu_{ij}, u^\nu_{kl}) = \delta_{\mu \nu} \delta_{ik} \delta_{jl}. 
$$

Both $ \CF^\infty(K_q) $ and $ \DF(K_q) $ admit universal $ C^* $-completions which we will denote by $ C(K_q) $ and $ C^*(K_q) $, respectively. 

The fundamental multiplicative unitary $ W $ is the algebraic multiplier of $ \CF^\infty(K_q) \otimes \DF(K_q) $ given by 
$$
W = \sum_{\mu \in \weights^+} \sum_{i,j = 1}^{\dim(V(\mu))} u^\mu_{ij} \otimes \omega^\mu_{ij}. 
$$
We may also view $ W $ as a unitary operator on the tensor product $ L^2(K_q) \otimes L^2(K_q) $, where $ L^2(K_q) $ is the GNS-Hilbert space 
of the Haar integral $ \phi $. 

Let us now define the Drinfeld double $ G_q = K_q \bowtie \hat{K}_q $. By definition, this is the algebraic quantum group given by the $ * $-algebra
$$
\CF^\infty_c(G_q) = \CF^\infty(K_q) \otimes \DF(K_q), 
$$
with comultiplication 
$$
\Delta_{G_q} = (\id \otimes \sigma \otimes \id)(\id \otimes \ad(W) \otimes \id)(\Delta \otimes \hat{\Delta}), 
$$
counit 
$$ 
\epsilon_{G_q} = \epsilon \otimes \hat{\epsilon}, 
$$
and antipode 
$$
S_{G_q}(f \otimes x) = W^{-1} (S(f) \otimes \hat{S}(x)) W = (S \otimes \hat{S})(W (f \otimes x) W^{-1}). 
$$
Here $ W $ denotes the multiplicative unitary from above. 
A positive left and right invariant Haar integral for $ \CF^\infty_c(G_q) $ is given by 
$$
\phi_{G_q}(f \otimes x) = \phi(f) \hat{\psi}(x), 
$$
as seen in \cite{PWlorentz}, \cite{VYcqg}. 

By taking the dual of $ \CF^\infty_c(G_q) $ in the sense of algebraic quantum groups, we obtain the convolution 
algebra $ \DF(G_q) = \DF(K_q) \bowtie \CF^\infty(K_q) $, which has $ \DF(K_q) \otimes \CF^\infty(K_q) $ 
as underlying vector space, equipped with the tensor product comultiplication and the multiplication 
\begin{align*}
(x \bowtie f)(y \bowtie g) &= x (y_{(1)}, f_{(1)}) y_{(2)} \bowtie f_{(2)} (\hat{S}(y_{(3)}), f_{(3)}) g.  
\end{align*}
The $ * $-structure of $ \DF(G_q) $ is defined in such a way that the natural inclusion homomorphisms $ \DF(K_q) \rightarrow \M(\DF(G_q)) $ and 
$ \CF^\infty(K_q) \rightarrow \M(\DF(G_q)) $ are $ * $-homomorphisms. Here $ \M(\DF(G_q)) $ denotes the algebraic multiplier algebra of $ \DF(G_q) $. 

Both algebras $ \CF^\infty_c(G_q) $ and $ \DF(G_q) $ admit universal $ C^* $-completions, which will 
be denoted by $ C_0(G_q) $ and $ C^*_\max(G_q) $, respectively. 

A unitary representation of $ G_q $ on a Hilbert space $ \H $ 
is defined to be a nondegenerate $ * $-homomorphism $ \pi: C^*_\max(G_q) \rightarrow \LH(\H) $. A basic example is the left regular representation 
of $ G_q $. It is obtained from a canonical $ * $-homomorphism $ C^*_\max(G) \rightarrow \LH(L^2(G_q)) $ 
where $ L^2(G_q) $ is the GNS-construction of the left Haar weight of $ G_q $. 
By definition, the reduced group $ C^* $-algebra $ C^*_\red(G_q) $ of $ G_q $ is the image of $ C^*_\max(G_q) $ under the left regular representation.

\section{The classical Baum-Connes deformation} \label{secbc}

In this section we review the definition of the Baum-Connes assembly map in the deformation picture. For background and more 
information we refer to section II.10 in \cite{Connesbook}, section 4 in \cite{BCH}, and \cite{Malowthesis}. 

Let us first recall a geometric construction called \emph{deformation to the normal cone} \cite{Connesbook}, \cite{Higsonmackeyanalogyktheory}, \cite{DSadiabatic},
which is a variant of a classical concept in algebraic geometry. We start with a smooth manifold $ M $ and a closed submanifold $ Z \subset M $, 
and write $ NZ $ for the normal bundle of this inclusion. The associated deformation to the normal cone is defined as 
$$
D_{Z \subset M} = \{0\} \times NZ \sqcup \mathbb{R}^\times \times M, 
$$
which becomes a smooth manifold in the following way. The topology and smooth structure are the obvious ones on
$ \{0\} \times NZ $ and $ \mathbb{R}^\times \times M $, and these subsets are glued together 
using an exponential map for the normal bundle. More precisely, assume that $ \exp: V \rightarrow W $ is a diffeomorphism such that $ \exp(z, 0) = z $ 
and $ d_{(z,0)} \exp = \id $ for all $ z \in Z $, where $ V \subset NZ $ is an open neighborhood of the zero 
section and $ W \subset M $ is an open set containing $ Z $. Then the map $ \theta $ given by 
$$
\theta(\tau, z, X) = 
\begin{cases}
(0, z, X) & \tau = 0 \\
(\tau, \exp(z, \tau X)) & \tau \neq 0 
\end{cases}
$$
is a diffeomorphism from a suitable neighborhood of $ \mathbb{R} \times Z \subset \mathbb{R} \times NZ $ onto 
$ \{0\} \times NZ \sqcup \mathbb{R}^\times \times W \subset D_{Z \subset M} $. 

In the sequel we will always restrict the parameter space in this construction from $ \mathbb{R} $ to $ [0,1] $, 
as this is better suited to the questions we are interested in. That is, we will consider $ \{0\} \times NZ \sqcup (0,1] \times M \subset D_{Z \subset M} $ and, by 
slight abuse of language, refer to this again as deformation to the normal cone. 
 
Now let $ G $ be a complex semisimple Lie group and let $ G = KAN $ be the Iwasawa decomposition of $ G $. 
We consider 
$$
\G_G = \{0\} \times K \times \mathfrak{k}^* \sqcup (0,1] \times G,  
$$
the deformation to the normal cone of the inclusion $ K \subset G $. Note here that the normal bundle $ NK $ is the trivial bundle 
$ K \times \mathfrak{an} $, and that $ \mathfrak{an} \cong \mathfrak{k}^* $ naturally. 
If we write $ G_0 = K \ltimes \mathfrak{k}^* $ for the semidirect product of $ K $ acting on the abelian Lie group $ \mathfrak{k}^* $ 
via the coadjoint action, then the fibrewise group structure on $ \G_G $ turns the latter into a Lie groupoid over the base space $ [0,1] $. 

The corresponding reduced groupoid $ C^* $-algebra $ C^*_\red(\G_G) $ is a continuous field of $ C^* $-algebras over $ [0,1] $ with 
fibres $ C^*(G_0) $ at $ 0 $ and $ C^*_\red(G) $ otherwise. 
Moreover, the bundle $ \G_G $ is trivial outside $ 0 $, and therefore induces a map $ K_*(C^*(G_0)) \rightarrow K_*(C^*_\red(G)) $ in $ K $-theory. 
Note that $ G_0 $ is amenable and that we have $ C^*(G_0) \cong K \ltimes C_0(\mathfrak{k}) $ by Fourier transform. 

The following result describes the deformation picture of the Baum-Connes assembly map, see section II.10 in \cite{Connesbook}.  

\begin{theorem} \label{CHdeformation}
Let $ G $ be a complex semisimple Lie group. There exists an isomorphism $ K^{\top}_*(G) \cong K_*(K \ltimes C_0(\mathfrak{k})) $ 
such that the diagram 
\begin{align*}
\xymatrix{
& K_*(K \ltimes C_0(\mathfrak{k})) \ar[dr] \\
K^{\top}_*(G) \ar[ur]_{\cong} \ar[rr]_{\mu} & & K_*(C^*_\red(G))
}
\end{align*}
is commutative, where $ \mu: K^{\top}_*(G) \rightarrow K_*(C^*_\red(G)) $ is the Baum-Connes assembly map
and the right downward arrow is the map arising from the deformation groupoid. 
\end{theorem} 

This result allows us to reinterpret the Baum-Connes assembly map in terms of the groupoid $ \G_G $. 
For a detailed account of its proof in the setting of arbitrary almost connected groups we refer to \cite{Malowthesis}.

\section{Quantization of the Cartan motion group} \label{seccmg}

In this section we explain how the classical Cartan motion group of $ G $ can be deformed into a quantum group. The corresponding 
continuous field of $ C^* $-algebras is a variant of the deformation to the normal cone construction. 

\subsection{Representations of crossed products by compact groups} 

Let us review some general facts regarding crossed product $ C^* $-algebras by actions of compact groups. For more information 
see \cite{RaeburnWilliamsmorita}, \cite{Williamscrossedproducts}. 

Let $ K $ be a compact group and let $ H \subset K $ be a closed subgroup. If $ V $ is a unitary representation of $ H $ 
then the induced representation is the Hilbert space 
$$
\ind_H^K(V) = L^2(K, V)^H 
$$
equipped with the left regular representation of $ K $, where $ H $ acts by $ (r \cdot \xi)(t) = r \cdot \xi(tr) $ on $ \xi \in L^2(K, V) $. 
Together with the action of $ C(K/H) $ on $ \ind_H^K(V) $ by pointwise multiplication, we obtain a covariant representation on $\ind_H^K(V)$. Hence $ \ind_H^K(V) $ is naturally a representation of $ K \ltimes C(K/H) $. 
In fact, the passage from $ V $ to $ \ind_H^K(V) $ corresponds to Rieffel induction under the Morita equivalence between $ C^*(H) $ 
and $ K \ltimes C(K/H) $.

Now let $ X $ be a locally compact space equipped with a continuous action of $ K $. 
It is well-known that the irreducible representations of the crossed product $ K \ltimes C_0(X) $ can be described as follows. 

\begin{prop} \label{crossedproductirreps}
Let $ K $ be a compact group acting on a locally compact space $ X $. Then all irreducible $ * $-representations of $ K \ltimes C_0(X) $ 
are of the form
$$
\ind_{K_x}^K(V)
$$
where $ x \in X $ and $ V $ is an irreducible unitary representation of the stabilizer group $ K_x $ of $ x $. Moreover, two representations $ \ind_{K_x}^K(V) $ and $ \ind_{K_y}^K(W) $ of this form are equivalent if and only if $ x = t \cdot y $ for some $ t \in K $ 
and $ W $ corresponds to $ V $ under the identification of $ K_x $ and $ K_y $ induced by conjugation with $ t $. 
\end{prop} 

\proof Note that $ \ind_{K_x}^K(V) $ is a representation of $ K \ltimes C_0(X) $ by first 
projecting $ K \ltimes C_0(X) \rightarrow K \ltimes C(K \cdot x) \cong K \ltimes C(K/K_x) $ and 
then considering the Morita equivalence between $ K \ltimes C(K/K_x) $ and $ C^*(K_x) $. 

Since Morita equivalence preserves irreducibility, it follows that inducing irreducible representations of stabilizers produces irreducible representations of 
the crossed product. In fact, all irreducible representations are obtained in this way, see for instance Proposition 8.7 in \cite{Williamscrossedproducts}. 

If $ y = t \cdot x $ for some $t \in K$, and $W$ corresponds to $V$ under $ \ad_t: K_x \rightarrow K_y, \ad_t(r) = trt^{-1} $, the corresponding isomorphism of stabilizer groups, 
then we obtain a unitary intertwiner $ T: \ind_{K_x}^K(V) \rightarrow \ind_{K_y}^K(W) $ by setting $ T(\xi)(x) = \xi(x t) $. 

Now assume that $ \ind_{K_x}^K(V) $ and $ \ind_{K_y}^K(W) $ are equivalent. If $ x $ and $ y $ are on the same $ K $-orbit, we can identify the stabilizers $K_x$ and $K_y$ via conjugation, and so we may assume $ x = y $. In this case we obtain $ V \cong W $ via the Morita equivalence between $ K \ltimes C(K/K_x) $ and $ C^*(K_x) $.

If $ x, y \in X $ are not on the same $ K $-orbit then the primitive ideals of $ K \ltimes C_0(X) $ corresponding to induced representations at $ x $ and $ y $ 
differ, and in particular there can be no nonzero intertwiner between such representations. \qed 

The following well-known fact regarding the structure of crossed products for actions of compact groups follows from the Imai-Takai biduality 
theorem \cite{ImaiTakaiduality}.

\begin{prop} \label{ttduality}
Let $ K $ be a compact group acting on a locally compact space $ X $. Then the crossed product $ K \ltimes C_0(X) $ 
is naturally isomorphic to the $ C^* $-algebra 
$$
C_0(X, \KH(L^2(K)))^K,  
$$
where $ K $ acts on $ C_0(X) $ with the given action and on $ \KH(L^2(K)) $ with the conjugation action with respect to the right regular 
representation on $ L^2(K) $. 
\end{prop} 

We note that under the identification in Proposition \ref{ttduality}, the (reduced) crossed product $ K \ltimes C_0(X) $ is 
the closed linear span of all elements of the form $ (y \otimes 1)\gamma(f) $ inside $ \KH(L^2(K)) \otimes C_0(X) $, where $ y \in C^*(K) $ and $ f \in C_0(X) $. 
Here $ \gamma: C_0(X) \rightarrow C(K) \otimes C_0(X) $ is given by $ \gamma(f)(t, x) = f(t \cdot x) $.

\subsection{The representation theory of $ K \ltimes C_0(\mathfrak{k}) $} 

Using the general method outlined above we shall now compute the irreducible representations of the Cartan motion group $ K \ltimes \mathfrak{k}^* $. 
For this it will be convenient to identify $ C^*(K \ltimes \mathfrak{k}^*) $ with the crossed product $ K \ltimes C_0(\mathfrak{k}) $. 
Unless explicitly stated otherwise, we assume throughout that $ K $ is a simply connected semisimple compact Lie group. 

We note first that each point of $ \mathfrak{k} $ is conjugate to a point of $ \mathfrak{t} $, see for instance
chapter IX in \cite{Bourbakilie7to9}. Therefore by Proposition \ref{crossedproductirreps} all irreducible representations of $ K \ltimes C_0(\mathfrak{k}) $ 
are of the form
\begin{displaymath}
\ind_{K_{X}}^K(V),
\end{displaymath}
where $ X \in \mathfrak{t} $ and $ V $ is an irreducible representation of the stabilizer subgroup $ K_X \subset K $ of $ X $. 
More concretely, the induced representation $ \ind_{K_X}^K(V) $ corresponds to the covariant pair given by 
the left regular representation of $ K $ on $ \ind_{K_X}^K(V) $ and 
\begin{displaymath}
\pi_{(V, X)}: C_0(\mathfrak{k}) \rightarrow \LH (\ind_{K_X}^K(V)), \quad \pi_{(V, X)}(f)( \xi)(k) = f(k \cdot X) \xi(k),
\end{displaymath}
where we write $ k \cdot X $ for the adjoint action of $ k $ on $ X $. 

If $ W $ denotes the Weyl group of $ K $, 
then the intersections of the orbits of $ K $ in $ \mathfrak{k} $ with $ \mathfrak{t} $ are the orbits of $ W $ in $ \mathfrak{t} $. 
For each $ X \in \mathfrak{t} $, the stabilizer subgroup $ K_X $ is a connected subgroup, 
and we have $ T \subset K_X \subset K $. 
Since the maximal torus of $ K_X $ is again $ T $ we have $ \Irr(K_X) = \Irr(T)/W_X $. 
Here $ W_X $ is the Weyl group of $ K_X $, and we note that this group identifies with the stabilizer of $ X $ in $ W $. 
Putting all these facts together we see that the spectrum of $ K \ltimes C_0(\mathfrak{k}) $, as a set, is
\begin{displaymath}
\bigsqcup_{X \in \mathfrak{t}/W} \weights/W_X \cong (\weights \times \mathfrak{t})/W. 
\end{displaymath}
Here we have identified $ \Irr(T) $ with the weight lattice $ \weights $.

In order to describe the algebra $ K \ltimes C_0(\mathfrak{k}) $ more concretely we use Proposition \ref{ttduality} to write
$$
K \ltimes C_0(\mathfrak{k}) \cong C_0(\mathfrak{k}, \KH(L^2(K)))^K. 
$$
Let $ \KH(L^2(K))^{K_X} $ be the fixed point subalgebra with respect to the right regular action of $ K_X $ on $ L^2(K) $. We note that 
if $ f \in C_0(\mathfrak{k}, \KH(L^2(K)))^K $, then $ f(X) \in \KH(L^2(K))^{K_X} $. Using the Peter-Weyl theorem one obtains 
\begin{displaymath}
\KH (L^2(K))^{K_X} \cong \bigoplus_{V \in \Irr(K_X)} \KH (L^2(K, V)^{K_X})
\end{displaymath}
by inducing the left regular representation of $ K_X $ to $ K $ and then applying Schur's Lemma.
In particular, if $ f \in C_0(\mathfrak{k}, \KH(L^2(K)))^K, X \in \mathfrak{t} $ and $ V $ is an irreducible representation of $ K_X $, 
then the projection of $ f(X) $ to the direct summand $ \KH(L^2(K, V)^{K_X}) $ in $\KH(L^2(K))^{K_X} $ corresponds to $ \ind_{K_X}^K(V) $. 

Consider the map 
$$ 
\pi: C_0(\mathfrak{k}, \KH(L^2(K)))^K \rightarrow C_0(\mathfrak{t}, \KH(L^2(K))) 
$$ 
obtained by restriction to $ \mathfrak{t} \subset \mathfrak{k} $. 
This map is injective since every orbit of the adjoint action of $ K $ on $ \mathfrak{k} $ meets $ \mathfrak{t} $. 

Our aim is to describe the image of $ \pi $. Given $ \mu \in \weights $ let us write $ V(\mu) $ for the corresponding representation of $ T $ and define
\begin{displaymath}
L^2(\E_{\mu}) = \ind_{T}^K(V(\mu)) = L^2(K, V(\mu))^T.
\end{displaymath}
The fixed point algebra in this case is
$$
\KH(L^2(K))^{T} = \bigoplus_{\mu \in \weights} \KH(L^2(\E_\mu)).
$$
Since $ T \subset K_X $ we have $ f(X) \in \KH(L^2(K))^{T} $ for all $ X \in \mathfrak{t} $. In particular, we see that $ \pi $ takes values 
in $ C_0(\weights \times \mathfrak{t}, \KH(\H))$, where $ \H = (\H_{\mu, X}) $ denotes the locally constant bundle of Hilbert spaces 
over $ \weights \times \mathfrak{t} $ with fibres $ \H_{\mu, X} = L^2(\E_\mu) $. Here $ C_0(\weights \times \mathfrak{t}, \KH(\H)) $ is 
the $ C^* $-algebra of sections of the $ C_0(\weights \times \mathfrak{t}) $-algebra of compact operators on the 
Hilbert $ C_0(\weights \times \mathfrak{t}) $-module corresponding to $ \H $.

The Weyl group $ W $ acts naturally on $ C_0(\weights \times \mathfrak{t}, \KH(\H)) = C_0(\mathfrak{t}, \bigoplus_{\mu \in \weights} \KH(L^2(\E_{\mu}))) $. 
Indeed, if $ w \in W $, there is a unitary isomorphism $ U_w: L^2(\E_{\mu}) \rightarrow L^2(\E_{w \mu}) $ 
defined by $ U_w(\xi)(k) = \xi(kw) $. 
Note that here we are viewing $ w \in K $, and strictly speaking we are picking a representative for an element of $ W $. However, the induced 
isomorphism $ \KH (L^2(\E_{\mu})) \cong \KH (L^2(\E_{w \mu})) $ is independent of this choice. The image of $ \pi $ is invariant 
under the action of the Weyl group, by the invariance under the action of $ K $. 

For $ X \in \mathfrak{t} $ and $ f \in C_0(\weights \times \mathfrak{t}, \KH(\H)) $ denote by $ f_X \in \KH(L^2(K)) $ 
the element obtained by evaluating $ f $ in the second variable and viewing the resulting section of the bundle $ \KH(\H) $ over $ \weights $ 
as an element of 
$$
\bigoplus_{\mu \in \weights} \KH(L^2(\E_\mu)) \subset \KH(L^2(K)). 
$$
With this notation in place we can describe the image of $ \pi $ as follows. 

\begin{theorem} \label{cmgplancherelform}
The $ C^* $-algebra $ K \ltimes C_0(\mathfrak{k}) $ is isomorphic to the subalgebra 
$$
A^L_0 = \{f \in C_0(\weights \times \mathfrak{t}, \KH(\H))^W \mid f_X \in \KH(L^2(K))^{K_X} \text{ for all } X \in \mathfrak{t} \} 
$$
of $ C_0(\weights \times \mathfrak{t}, \KH(\H))^W $. 
\end{theorem} 

\proof We already observed that the map $ \pi $ defined above is injective and takes values in $ A^L_0 $. For surjectivity, we start by understanding the form of the irreducible representations of $ A^L_0 $.

Let $ \pi: A^L_0 \rightarrow \LH (\V) $ be an irreducible representation of $ A^L_0 $ on a Hilbert space $ \V $. By standard $ C^* $-algebra theory, 
there exists an irreducible representation $ \Pi: C_0(\weights \times \mathfrak{t}, \KH (\H)) \rightarrow \LH (\K) $ on some Hilbert space $ \K $ and 
a $ A^L_0 $-invariant subspace $ \L \subset \K $ such that the restriction of $ \Pi $ to $ A^L_0 $, acting on $ \L $, is equivalent to $ \pi $. 
In particular, we can assume $ \V \subset \K $ such that the action of $ \pi $ on $ \V $ is given by $ \Pi $. Now note that 
$ C_0(\weights \times \mathfrak{t}, \KH (\H)) = C_0(\mathfrak{t}, \bigoplus_{\mu \in \weights} \KH (L^2(\E_{\mu}))) $. Then, up to isomorphism, 
$ \Pi $ must be given by point evaluation at an element $ X \in \mathfrak{t} $, followed by projection to a direct summand 
of $ \oplus_{\mu \in \weights} \KH (L^2(\E_{\mu})) $. 

In particular, each irreducible representation of $ A^L_0 $ factors through point evaluation at some element of $ \mathfrak{t} $. If $ f \in A^L_0 $, 
then $ f_X \in \KH(L^2(K))^{K_X} $, which we saw has a direct sum decomposition
\begin{displaymath}
\KH(L^2(K))^{K_X} \cong \bigoplus_{V \in \Irr(K_X)} \KH (L^2(K, V)^{K_X}).
\end{displaymath}
We conclude that any irreducible representation of $ A^L_0 $ is given by point evaluation at an element $ X \in \mathfrak{t} $, followed by projection to one 
of these direct summands of $ \KH(L^2(K))^{K_X} $. Notice that restricting an irreducible representation of $ A^L_0 $ to the image of $ \pi $ gives an 
irreducible representation of the latter, by our concrete description of the irreducible representations of $ C_0(\mathfrak{k}, \KH(L^2(K)))^K $. 

The result then follows from Dixmier's Stone-Weierstrass theorem, see section 11.1 in \cite{Dixmiercstar}. \qed

In order to illustrate Theorem \ref{cmgplancherelform} let us consider the special case $ K = SU(2) $. In this case the stabilizer group $ K_X $ is 
equal to $ T $ if $ X \in \mathfrak{t} \cong \mathbb{R} $ is non-zero, and $ K_X = K $ if $ X = 0 $. 
It follows that the algebra described in Theorem \ref{cmgplancherelform} identifies with 
$$
A^L_0 = \{f \in C_0(\mathbb{Z} \times \mathbb{R}, \KH(\H))^{\mathbb{Z}_2} \mid f_0 \in \bigoplus_{\nu \in \weights^+} \KH(V(\nu)) \subset \KH(L^2(K)) \}, 
$$
where $ \KH(V(\nu)) $ is represented via the left regular representation of $ K $.

\subsection{The representation theory of $ K \ltimes C(K) $} 

Let us next compute the irreducible representations of the $ C^* $-algebra $ K \ltimes C(K) $ of the quantum Cartan motion group. 
The arguments are parallel to the ones in the previous subsection. 
We remark that $ K \ltimes C(K) $ is the group $ C^* $-algebra of the Drinfeld double of the classical group $ K $. 

The irreducible representations of $ K \ltimes C(K) $ are labelled by $ (\weights \times T)/W $. Indeed, each point of $ K $ is conjugate to a 
point of $ T $, and so Proposition \ref{crossedproductirreps} tells us all irreducible representations of $ K \ltimes C(K) $ are of the form
\begin{displaymath}
\ind_{K_{t}}^K(V),
\end{displaymath}
where $ t \in T $ and $ V $ is an irreducible representation of the centralizer subgroup $ K_t \subset K $ of $ t $ in $ K $. 
The intersections of the conjugacy classes in $ K $ with $ T $ are the orbits of $ W $ on $ T $. 
For each $ t \in T $, the centralizer subgroup $ K_t $ is a connected subgroup, 
and we have $ T \subset K_t \subset K $. Therefore $ \Irr(K_t) = \Irr(T)/W_t $, where $ W_t $ is the Weyl group of $ K_t $, which identifies with the stabilizer 
of $ t $ in $ W $. 
We conclude that the spectrum of $ K \ltimes C(K) $, as a set, is
\begin{displaymath}
\bigsqcup_{t \in T/W} \weights/W_t \cong (\weights \times T)/W.
\end{displaymath}

In the same way as in the analysis of the classical Cartan motion group, let us consider the isomorphism 
$$
K \ltimes C(K) \cong C(K, \KH(L^2(K)))^K 
$$
obtained from Proposition \ref{ttduality}. If $ f \in C(K, \KH(L^2(K)))^K $ then $ f(t) $ 
is contained in $ \KH(L^2(K))^{K_t} $, and each irreducible representation of $ C(K, \KH(L^2(K)))^K $ 
can be described as point evaluation at some $ t \in T $ followed by projection onto a direct summand of 
$$ 
\KH(L^2(K))^{K_t} \cong \bigoplus_{V \in \Irr(K_t)} \KH (L^2(K, V)^{K_t}).
$$
The map
$$ 
\pi: C(K, \KH(L^2(K)))^K \rightarrow C(T, \KH(L^2(K)))
$$ 
obtained by restricting operator-valued functions to $ T \subset K $ is injective,  
and the image of $ \pi $ is contained in $ C_0(\weights \times T, \KH(\H))^W $, where $ \H = (\H_{\mu, t}) $ denotes the 
locally constant bundle of Hilbert spaces over $ \weights \times T $ with fibres $ \H_{\mu, t} = L^2(\E_\mu) $. 

For $ t \in T $ and $ f \in C_0(\weights \times T, \KH(\H))^W $ denote by $ f_t \in \KH(L^2(K)) $ the element obtained by evaluating $ f $ in the second variable and viewing the resulting section of the bundle $ \KH(\H) $ over $ \weights $ as an element of 
$$
\bigoplus_{\mu \in \weights} \KH(L^2(\E_\mu)) \subset \KH(L^2(K)). 
$$
We then have the following analogue of Theorem \ref{cmgplancherelform}, with essentially the same proof. 

\begin{theorem} \label{qcmgplancherelform}
The $ C^* $-algebra $ K \ltimes C(K) $ is isomorphic to the subalgebra 
$$
A^L_1 = \{f \in C_0(\weights \times T, \KH(\H))^W \mid f_t \in \KH(L^2(K))^{K_t} \text{ for all } t \in T \} 
$$
of $ C_0(\weights \times T, \KH(\H))^W $. 
\end{theorem} 

Let us again consider explicitly the case $ K = SU(2) $. In this case the centralizer group $ K_t $ is equal to $ K $ if $ t = \pm 1 $ 
and $ K_t = T \cong S^1 $ otherwise. 

For $ t \in T $ let us write 
$$ 
V_{\mu, t} = \ind_{K_t}^K(V(\mu))
$$
for the induced representation of the irreducible representation $ V(\mu) $ of $ K_t $ associated to the weight $ \mu \in \weights $. 
Here we tacitly declare $ V(\mu) = V(-\mu) $ if $ \mu \in \weights^+ $ and $ t = \pm 1 $. Two representations $ V_{\mu, t}, V_{\mu', t'} $ 
are equivalent if and only if $ (\mu', t') = (\pm \mu, t^{\pm 1}) $. 

Let us write down the representations $ V_{\mu, t} $ explicitly. 
If $ t \neq \pm 1 $ then $ V_{\mu, t} $ is the space $ L^2(\E_\mu) $ with the action of $ C(K) $ given by 
$$
(f \cdot \xi)(k) = f(k t k^{-1}) \xi(k), 
$$
and the action of $ K $ given by 
$$
(s \cdot \xi)(k) = \xi(s^{-1} k).
$$
If $ t = \pm 1 $ we have $ V_{\mu, t} = V(\mu) $ and the action of $ C(K) $ reduces to
$$
f \cdot \xi = f(t) \xi. 
$$
The algebra described in Theorem \ref{qcmgplancherelform} is  
$$
A^L_1 = \{f \in C_0(\mathbb{Z} \times S^1, \KH(\H))^{\mathbb{Z}_2} \mid f_{\pm 1} \in \bigoplus_{\nu \in \weights^+} \KH(V(\nu)) \subset \KH(L^2(K)) \}, 
$$
where $ \KH(V(\nu)) $ is again represented via the left regular representation of $ K $.

\subsection{The quantization field for the Cartan motion group} 

In this section we discuss the quantization of the classical Cartan motion group, represented by the $ C^* $-algebra $ K \ltimes C_0(\mathfrak{k}) $, 
towards the quantum Cartan motion group, which is represented by the $ C^* $-algebra $ K \ltimes C(K) $. 

Let us consider the deformation to the normal cone construction outlined in section \ref{secbc} for the manifold $ K $ and its one-point submanifold 
consisting of the identity element of $ K $. We obtain a Lie groupoid $ \G_K $ over the base space $ [0,1] $, also known as the corresponding adiabatic 
groupoid \cite{DSadiabatic}, with object space 
$$
\G_K = \{0\} \times \mathfrak{k} \sqcup (0,1] \times K, 
$$
equipped with the topology and smooth structure as defined in section \ref{secbc}.
The groupoid structure is given by the fibrewise group operations. 

For our purposes it is convenient to describe this construction directly on the level of functions. 
Given any $ f_0 \in C_c(\mathfrak{k}) $ we obtain a continuous function $ f $ on $ \G_K $ by considering the family $ (f_\tau)_{\tau \in [0,1]} $ of functions 
defined as follows. We set $ f_\tau(\exp(X)) = f(\tau^{-1} X) $ for $ \tau > 0 $ small and $ X $ in a neighbourhood of $ 0 $ 
on which the exponential map is a diffeomorphism. 
Extending by zero and cutting off in the $ \tau $-direction we obtain continuous functions $ f_\tau \in C(K) $ for all $ \tau > 0 $ in this way. 
Evaluating the family $ (f_\tau)_{\tau \in [0,1]} $ in the fibres of $ \G_K $ then yields $ f \in C_0(\G_K) $.  
We also consider the continuous functions on $ \G_K $ obtained from all families $ (f_\tau)_{\tau \in [0,1]} $ such that $ f_\tau = 0 $ for $ \tau $ small 
and $ (\tau, x) \mapsto f_\tau(x) $ is continuous on $ (0,1] \times K $. 

By definition of the topology of $ \G_K $, the collection of all functions described above yields a $ * $-algebra which is dense in the 
$ C^* $-algebra $ C_0(\G_K) $ of functions on $ \G_K $. This shows in particular that $ C_0(\G_K) $ is a continuous field of $ C^* $-algebras over $ [0,1] $
with fibres 
$$
C_0(\G_K)_\tau = 
\begin{cases} 
C_0(\mathfrak{k}) & \text{ if } \tau = 0 \\
C(K) & \text{ else.}
\end{cases}
$$
Using the fact that continuous fields are preserved by taking crossed products with respect to fibrewise actions of amenable groups, 
see \cite{PRtwistedcp2}, we obtain the following result. 

\begin{prop} \label{cmgdeformation}
Let $ K $ be a connected compact Lie group and let $ \G_K $ be the associated adiabatic groupoid as above. 
The crossed product $ L = K \ltimes C_0(\G_K) $ with respect to the fibrewise adjoint action 
is a continuous field of $ C^* $-algebras with fibres $ L_0 = K \ltimes C_0(\mathfrak{k}) $ and $ L_\tau = K \ltimes C(K) $ for $ \tau > 0 $. 
\end{prop} 

We call the continuous field constructed in Proposition \ref{cmgdeformation} the \emph{quantization field for the Cartan motion group}. 

Using Theorem \ref{cmgplancherelform} and Theorem \ref{qcmgplancherelform} let us define another continuous field of $ C^* $-algebras $ A^L $ over $ [0,1] $ 
whose fibres are isomorphic to the fibres of the quantization field $ L $. More specifically, we stipulate that the field $ A^L $ is trivial away from $ 0 $, and if 
$ f_0 \in K \ltimes C_0(\mathfrak{k}) \subset C_0(\weights \times \mathfrak{t}, \KH(\H))^W $ 
has compact support as a function on $ \weights \times \mathfrak{t} $, 
then a continuous section $ (f_\tau)_{\tau \in [0,1]} $ of $ A^L $ is obtained by setting 
$$
f_\tau(\mu, \exp(X)) = f(\mu, \tau^{-1} X)
$$
for $ \tau > 0 $ small and $ X $ in a neighbourhood of $ 0 $ on which $ \exp $ is a diffeomorphism, extended by $ 0 $ to all of $ T $. 
Note here that $ f_\tau $ for $ \tau > 0 $ defines an element of $ K \ltimes C(K) \subset C_0(\weights \times T, \KH(\H))^W $. 
The fibre at $ \tau = 0 $ of the field $ A^L $ identifies with $ A^L_0 $ as in Theorem \ref{cmgplancherelform}, and the fibres of 
$ A^L $ at points $ \tau > 0 $ can be identified with $ A^L_1 $ as in Theorem \ref{qcmgplancherelform}.  

\begin{prop} 
The continuous field $ A^L $ described above is isomorphic to the quantization field of the Cartan motion group. 
\end{prop} 

\proof Let $ f $ be an element of the convolution algebra $ C(K, C_0(\G_K)) \subset K \ltimes C_0(\G_K) $ of the form $ f = g \otimes h $ 
where $ g \in C(K) $ and $ h \in C_0(\G_K) $ is the generating section associated with a compactly supported
function $ h_0 $ in $ C_c(\mathfrak{k}) $ as in the description of $ C_0(\G_K) $ explained further above. 
Since linear combinations of elements of this type form a dense subspace of $ L $ it suffices to show that $ f $ defines a continuous section of $ A^L $. 

Under the isomorphism provided by Theorem \ref{cmgplancherelform}, we have
\begin{displaymath}
((g \otimes h_0)(\mu, X))(\xi)(k) =  \int_{K} g(s) h(s^{-1}k \cdot X) \xi(s^{-1}k) ds
\end{displaymath}
for $ \mu \in \weights, X \in \mathfrak{t}, \xi \in L^2(\E_{\mu}) $ and $ k \in K $. Similarly, under the isomorphism provided by 
Theorem \ref{qcmgplancherelform}, we have
\begin{displaymath}
((g \otimes h_{\tau})(\mu, \exp(X)))(\xi)(k) =  \int_{K} g(s)h(\tau^{-1} s^{-1}k \cdot X) \xi(s^{-1}k) ds 
\end{displaymath}
for $ \tau > 0 $ small. From these formulae we can see that under these isomorphisms on each fibre, the element $ f $ is mapped to a continuous 
section of $ A^L $, as required. \qed

\section{Quantization of complex semisimple Lie groups} \label{seccg}

In this section we describe a deformation of the reduced group $ C^* $-algebra of a complex semisimple Lie group $ G $. 
The resulting continuous field of $ C^* $-algebras is another instance of the deformation to the normal cone construction. 

\subsection{The tempered representations of classical complex groups} 

In this subsection we review the structure of the space of irreducible tempered representations of complex semisimple groups. 
Throughout we fix a simply connected complex semisimple Lie group $ G $ with maximal compact subgroup $ K $, and a compatible choice of Borel 
subgroup $ B \subset G $ and maximal torus $ T \subset K $. 

Let $ \mu \in \weights $. Then the space of smooth sections $ \Gamma(\E_\mu) \subset L^2(\E_\mu) $ of the induced vector bundle 
$ \E_\mu = K \times_T \mathbb{C}_\mu $ over $ G/B = K/T $ corresponding to $ \mu $ is the subspace of $ C^\infty(K) $ 
of weight $ \mu $ with respect to the right translation action of $ T $ given by 
$$
(t \hit \xi)(x) = \xi(xt). 
$$
Equivalently, we have 
$$
\Gamma(\E_\mu) = \{\xi \in C^\infty(K) \mid (\id \otimes \pi_T) \Delta(\xi) = \xi \otimes z^\mu \}, 
$$ 
where $ \pi_T: C^\infty(K) \rightarrow C^\infty(T) $ is the projection homomorphism and $ z^\mu \in C^\infty(T) $ 
is the unitary corresponding to the weight $ \mu $.  

Given $ \lambda \in \mathfrak{t}^* $ we can identify 
$$ 
\Gamma(\E_\mu) = \{\xi \in C^\infty(G) \mid (\id \otimes \pi_B) \Delta(\xi) = \xi \otimes z^\mu \otimes \chi_{\lambda + 2 \rho} \} \subset C^\infty(G),  
$$
where we view $ G = KAN = K \times \exp(\mathfrak{a}) \times \exp(\mathfrak{n}) $ and write $ \chi_{\lambda + 2 \rho} $ for the character 
on $ AN $ given by 
$$ 
\chi_{\lambda + 2 \rho}(\exp(a) \exp(n)) = e^{(\lambda + 2 \rho, a)}; 
$$
recall that $ \rho $ denotes the half sum of all positive roots. 
The left regular action of $ G $ on $ \Gamma(\E_\mu) = \Gamma(\E_{\mu, \lambda}) $ defines a representation of $ G $ which is called 
the \emph{principal series representation} with parameter $ (\mu, \lambda) \in \weights \times \mathfrak{t}^* $.  
Since $ \lambda \in \mathfrak{t}^* $ the representation $ \Gamma(\E_{\mu, \lambda}) \subset C^\infty(K) $ is unitary for the standard scalar 
product on $ C^\infty(K) $ induced from the Haar measure on $ K $. In particular, we obtain a corresponding 
nondegenerate $ * $-representation $ \pi_{\mu, \lambda}: C^*_\max(G) \rightarrow \LH(\H_{\mu, \lambda}) $, 
where $ \H_{\mu, \lambda} \subset L^2(K) $ is the Hilbert space completion of $ \Gamma(\E_{\mu, \lambda}) $. 

The following result is due to \v Zelobenko \cite{Zelobenkoirreducibility} and Wallach \cite{Wallachcyclicvectors}. 

\begin{theorem} 
For all $ (\mu, \lambda) \in \weights \times \mathfrak{t}^* $ the unitary principal series representation $ \H_{\mu, \lambda} $ is an 
irreducible unitary representation of $ G $. 
\end{theorem} 

The Weyl group $ W $ acts on the parameter space $ \weights \times \mathfrak{t}^* $ by 
$$
w (\mu, \lambda) = (w \mu, w \lambda). 
$$
Using this action one can describe the isomorphisms between unitary principal series representations as follows, see \cite{KSuniformlyboundedrepIII}. 

\begin{theorem} \label{thmweylgroupactionirreds}
Let $ (\mu, \lambda), (\mu', \lambda') \in \weights \times \mathfrak{t}^* $. Then $ \H_{\mu, \lambda} $ 
and $ \H_{\mu', \lambda'} $ are equivalent representations of $ G $ if and only if $ (\mu, \lambda), (\mu', \lambda') $ 
are on the same Weyl group orbit, that is, if and only if 
$$ 
(\mu', \lambda') = (w \mu, w \lambda) 
$$ 
for some $ w \in W $. 
\end{theorem} 

Combining these results with Harish-Chandra's Plancherel theorem \cite{HarishChandraplancherelcomplex} one obtains the following description of 
the reduced group $ C^* $-algebra of $ G $. 

\begin{theorem} \label{reducedcstarstructure}
Let $ G $ be a complex semisimple Lie group, and let $ \H = (\H_{\mu, \lambda})_{\mu, \lambda} $ be the 
Hilbert space bundle of unitary principal series representations of $ G $ over $ \weights \times \mathfrak{t}^* $. Then one obtains an isomorphism 
$$
C^*_\red(G) \cong C_0(\weights \times \mathfrak{t}^*, \KH(\H))^W 
$$
induced by the canonical $ * $-homomorphism $ \pi: C^*_\max(G) \rightarrow C_0(\weights \times \mathfrak{t}^*, \KH(\H)) $. 
\end{theorem}

\subsection{The tempered representations of complex quantum groups} 

In this subsection we review some aspects of the representation theory of complex quantum groups. For more details and background 
we refer to \cite{VYcqg}. 

Throughout we fix $ q = e^h $ such that $ 0 < q < 1 $ and set 
$$
\mathfrak{t}^*_q = \mathfrak{t}^*/i \hbar^{-1} \roots^\vee,
$$
where $ \hbar = \frac{h}{2 \pi} $ and $ \roots^\vee $ is the coroot lattice of $ \mathfrak{g} $. 

Let $ \mu \in \weights $. Then we define the space of sections $ \Gamma(\E_\mu) \subset \CF^\infty(K_q) $ of the induced vector bundle 
$ \E_\mu $ corresponding to $ \mu $ to be the subspace of $ \CF^\infty(K_q) $ of weight $ \mu $ with respect to the $ \DF(K_q) $-module 
structure 
$$
x \hit \xi = \xi_{(1)} (x, \xi_{(2)}).  
$$
Equivalently, we have 
$$
\Gamma(\E_\mu) = \{\xi \in \CF^\infty(K_q) \mid (\id \otimes \pi_T) \Delta(\xi) = \xi \otimes z^\mu \}, 
$$ 
where $ \pi_T: \CF^\infty(K_q) \rightarrow \CF^\infty(T) $ is the canonical projection homomorphism and $ z^\mu \in \CF^\infty(T) = \mathbb{C}[\weights] $ 
is the generator corresponding to the weight $ \mu $.  

Let $ \mathfrak{t}_q^* = \mathfrak{t}^*/i\hbar^{-1} \roots^\vee $ where $ \roots^\vee $ is the coroot lattice. 
For $ \lambda \in \mathfrak{t}_q^* $ we define the twisted left adjoint representation of $ \CF^\infty(K_q) $ on $ \Gamma(\E_\mu) $ by 
\begin{equation*} 
f \cdot \xi = f_{(1)} \, \xi \, S(f_{(3)}) (K_{2 \rho + \lambda}, f_{(2)}).
\end{equation*}
Together with the left coaction $ \Gamma(\E_\mu) \rightarrow \CF^\infty(K_q) \otimes \Gamma(\E_\mu) $ given by comultiplication,  
this turns $ \Gamma(\E_\mu) $ into a Yetter-Drinfeld module. 

It is convenient to switch from the left coaction on $ \Gamma(\E_\mu) $ to the left $ \DF(K_q) $-module structure given by 
$$
x \cdot \xi = (\hat{S}(x), \xi_{(1)}) \xi_{(2)} 
$$
for $ x \in \DF(K_q) $.   
The space $ \Gamma(\E_\mu) = \Gamma(\E_{\mu, \lambda}) $ with the actions as above is called 
the \emph{principal series Yetter-Drinfeld module} with parameter $ (\mu, \lambda) \in \weights \times \mathfrak{t}^*_q $.  
Since $ \lambda \in \mathfrak{t}^*_q $ the Yetter-Drinfeld module $ \Gamma(\E_{\mu, \lambda}) \subset \CF^\infty(K_q) $ is unitary for the scalar 
product induced from the Haar state on $ \CF^\infty(K_q) $. In particular, we obtain a corresponding 
nondegenerate $ * $-representation $ \pi_{\mu, \lambda}: C^*_\max(G_q) \rightarrow \LH(\H_{\mu, \lambda}) $, 
where $ \H_{\mu, \lambda} \subset L^2(K_q) $ is the Hilbert space completion of $ \Gamma(\E_{\mu, \lambda}) $. 

The unitary representations of $ G_q $ on $ \H_{\mu, \lambda} $ for $ (\mu, \lambda) \in \weights \times \mathfrak{t}^*_q $ as above 
are called unitary principal series representations. 

For a proof of the following result we refer to \cite{VYcqg}. 

\begin{theorem} 
For all $ (\mu, \lambda) \in \weights \times \mathfrak{t}^*_q $ the unitary principal series representation $ \H_{\mu, \lambda} $ is an 
irreducible unitary representation of $ G_q $. 
\end{theorem} 

As in the classical situation, there are nontrivial intertwiners between unitary principal series representations of $ G_q $.  
The classical Weyl group $ W $ acts on the parameter space $ \weights \times \mathfrak{t}^*_q $ by 
$$
w (\mu, \lambda) = (w \mu, w \lambda). 
$$
The following result describes the isomorphisms between unitary principal series representations in the quantum case \cite{VYcqg}. 

\begin{theorem} \label{thmweylgroupactionirredsquantum}
Let $ (\mu, \lambda), (\mu', \lambda') \in \weights \times \mathfrak{t}^*_q $. Then $ \H_{\mu, \lambda} $ 
and $ \H_{\mu', \lambda'} $ are equivalent representations of $ G_q $ if and only if $ (\mu, \lambda), (\mu', \lambda') $ 
are on the same Weyl group orbit, that is, if and only if 
$$ 
(\mu', \lambda') = (w \mu, w \lambda) 
$$ 
for some $ w \in W $. 
\end{theorem} 

Let $ \H = (\H_{\mu, \lambda})_{\mu, \lambda} $ be the Hilbert space bundle of unitary principal series representations of $ G_q $ 
over $ \weights \times \mathfrak{t}^*_q $. 
Then we obtain a $ * $-homomorphism $ \pi: C^*_\max(G) \rightarrow C_0(\weights \times \mathfrak{t}^*_q, \KH(\H)) $ 
by setting $ \pi(x)(\mu, \lambda) = \pi_{\mu, \lambda}(x) $, and this map is used to describe the structure of 
the reduced $ C^* $-algebra of $ G_q $ as follows \cite{VYcqgplancherel}. 

\begin{theorem} \label{reducedcstarstructurequantum}
Let $ G_q $ be a complex semisimple quantum group, and let $ \H = (\H_{\mu, \lambda})_{\mu, \lambda} $ be the 
Hilbert space bundle of unitary principal series representations of $ G_q $ over $ \weights \times \mathfrak{t}^*_q $. Then one obtains an isomorphism 
$$
C^*_\red(G_q) \cong C_0(\weights \times \mathfrak{t}^*_q, \KH(\H))^W 
$$
induced by the canonical $ * $-homomorphism $ \pi: C^*_\max(G_q) \rightarrow C_0(\weights \times \mathfrak{t}^*_q, \KH(\H)) $. 
\end{theorem}

In the case $ G_q = SL_q(2, \mathbb{C}) $, the statement of Theorem \ref{reducedcstarstructurequantum} follows already from 
the work of Pusz-Woronowicz \cite{Puszunitaryrepresentations}, \cite{PWquantumlorentzgelfand} and Buffenoir-Roche \cite{BRLorentz}.

\subsection{The quantization field of a complex group} 

Using the structure of the reduced group $ C^* $-algebras $ C^*_\red(G) $ and $ C^*_\red(G_q) $ described above we shall now 
construct a continuous field of $ C^* $-algebras relating them. 

Consider the parameter space $ M = \weights \times \mathfrak{t}^*_q $ of unitary principal series representations of the quantum group $ G_q $. 
If we let $ Z = \weights \times \{0\} \subset M $ then we can view the normal bundle $ NZ = \weights \times \mathfrak{t}^* $ as the parameter 
space of unitary principal series representations of the classical group $ G $. 
Moreover, due to Peter-Weyl theory the underlying Hilbert spaces of classical and quantum unitary principal series representations associated with 
the same parameter $ \mu \in \weights $ can be naturally identified. Using the same techniques as in section \ref{seccmg} we obtain 
a continuous field of $ C^* $-algebras $ B $ over $ [0,1] $ with fibres $ B_0 = C_0(\weights \times \mathfrak{t}^*, \KH(\H)) $ 
and $ B_\tau = C_0(\weights \times \mathfrak{t}^*/(i \hbar^{-1}\tau^{-1} \roots^\vee), \KH(\H)) $ for $ \tau > 0 $. 

The Weyl group action on the bundle $ \H $ arising from Theorem \ref{thmweylgroupactionirreds} and Theorem \ref{thmweylgroupactionirredsquantum} 
is compatible with this continuous field structure, see the explicit formulas for intertwining operators obtained in section 5 of \cite{VYcqg}. 
Therefore we obtain an action of $ W $ on $ B $. 
Taking the $ W $-invariant part of $ B $ yields a continuous field of $ C^* $-algebras $ R = B^W $ over $ [0,1] $ with 
fibres $ R_0 = C^*_\red(G) $ and $ R_\tau = C^*_\red(G_{q^\tau}) $ for $ \tau > 0 $. We shall refer to this 
continuous field as the \emph{quantization field} of $ G $. 

For our purposes it will be convenient to reparameterize the spaces of principal series representations of $ G $ and $ G_{q^\tau} $, respectively.
Let $ X_*(T) = \ker(\exp: \mathfrak{t} \rightarrow T) $ denote the coweight lattice of $ T $. 
Then we obtain a linear isomorphism $ \gamma: \mathfrak{t} \rightarrow \mathfrak{t}^* $ 
determined by $ (\gamma(X), \mu) = \mu(X) $ for all $ \mu \in \weights $. 
Similarly, let us write $ \mathfrak{t}_{q^\tau} = \mathfrak{t}/h^{-1} \tau^{-1} X_*(T) $. Using this notation we obtain a linear 
isomorphism $ \gamma_{q^\tau}: \mathfrak{t}_{q^\tau} \rightarrow \mathfrak{t}^*_{q^\tau} $ 
determined by the same formula as $ \gamma $. 

These isomorphisms are $ W $-equivariant.  
According to Theorem \ref{reducedcstarstructure} and Theorem \ref{reducedcstarstructurequantum}, we can therefore identify the fibres of 
the quantization field of $ G $ with $ A^R_0 = C_0(\weights \times \mathfrak{t}, \KH(\H))^W $ 
and $ A^R_\tau = C_0(\weights \times T, \KH(\H))^W $ for $ \tau > 0 $, taking into account the 
canonical rescaling $ \mathfrak{t}_{q^\tau} = \mathfrak{t}/h^{-1} \tau^{-1} X_*(T) \cong \mathfrak{t}/X_*(T) \cong T $ 
for all $ \tau > 0 $. Let us write $ A^R $ for the resulting continuous field. 

The quantization field of $ G $ can be viewed as a noncommutative deformation to the normal cone construction applied to the tempered dual of $ G_q $.

\section{The quantum Baum-Connes assembly map} \label{secqbc}

In this section we define a quantum analogue of the Baum-Connes assembly field and show that it induces an isomorphism in $ K $-theory. 

\subsection{The quantum Baum-Connes field} 

The key ingredient in our approach is to obtain a continuous field of $ C^* $-algebras by varying the deformation parameter 
in the construction of the Drinfeld double. More precisely, we shall fix $ q = e^h $ and construct a continuous field $ Q = C^*_\red({\bf G}) $ 
of $ C^* $-algebras over $ [0,1] $ with fibres $ Q_\sigma = C^*_\red(G_{q^\sigma}) $. 

Let us first recall that the family of $ C^* $-algebras $ (C(K_{q^\sigma}))_{\sigma \in [0,1]} $ assembles into a continuous field $ C({\bf K}) $ 
of $ C^* $-algebras over $ [0,1] $ as follows \cite{NTKhomologydirac}. If $ V(\mu) $ denotes the underlying Hilbert space of the 
irreducible representation of $ K $ of highest weight $ \mu \in \weights^+ $, and $ V(\mu)_\sigma $ the corresponding irreducible representation 
of $ K_{q^\sigma} $, then we may fix a continuous family of unitary isomorphisms $ V(\mu)_\sigma \cong V(\mu) $ 
which is the identity for $ \sigma = 0 $. Here continuity means that the unbounded multipliers $ E^{q^\sigma}_i, F^{q^\sigma}_i, H^{q^\sigma}_i $ 
of $ C^*(K_{q^\sigma}) $ define continuous families of operators in $ \KH(V(\mu)) $ for $ \sigma \in [0,1] $ under these isomorphisms. 
 
Then for each $ v \in V(\mu)^*, w \in V(\mu) $ the matrix element $ \bra v| \bullet| w \ket $ naturally becomes a continuous section of the 
field $ C({\bf K}) $, and this in turn determines the continuous field structure of $ C({\bf K}) $. 

Dually, note that the $ C^* $-algebras $ C^*(K_{q^\sigma}) $ naturally identify with $ C^*(K) = \bigoplus_{\mu \in \weights^+} \KH(V(\mu)) $
under the above continuous family of isomorphisms. In particular, these algebras assemble into a constant continuous field $ C^*({\bf K}) $ 
over $ [0,1] $ in an obvious way. 
 
Roughly, in order to construct the field $ Q = C^*_\red({\bf G}) $ we need to combine the fields $ C({\bf K}) $ and $ C^*({\bf K}) $ as in the Drinfeld double 
construction. 

More precisely, consider the $ C^* $-algebra $ E = \prod_{\sigma \in [0,1]} C^*_\red(G_{q^\sigma}) $ consisting of all uniformly norm-bounded families of elements 
in $ C^*_\red(G_{q^\sigma}) $. Let $ C^*_\red({\bf G}) \subset E $ be the $ C^* $-algebra generated by 
all operators of the form $ (\omega^\mu_{ij} \bowtie u^\nu_{kl}) f $ where $ f \in C[0,1], \mu, \nu \in \weights^+ $, and the indices run over all allowed values. 
Note that $ C^*_\red({\bf G}) $ is indeed contained in $ E $ since the norms of $ \omega^\mu_{ij} $ and $ u^\nu_{kl} $ 
are uniformly bounded by $ 1 $; in the case of $ \omega^\mu_{ij} $ this holds because the elements $ \omega^\mu_{ij} $ for fixed $ \mu $ form a 
full matrix algebra, and for $ u^\nu_{kl} $ it suffices to note that the matrix $ u^\nu = (u^\nu_{kl}) $ over $ C(K_{q^\sigma}) $ is unitary 
for all $ \sigma $. 

By construction, the $ C^* $-algebra $ C^*_\red({\bf G}) $ is a $ C[0,1] $-algebra. 
In order to describe this algebra further it is convenient to reparameterize the space of principal series representations of $ G_{q^\sigma} $ as 
described towards the end of section \ref{seccg}. 
Then according to Theorem \ref{qcmgplancherelform} and Theorem \ref{reducedcstarstructurequantum} we obtain a natural inclusion 
$$ 
C^*_\red({\bf G}) \subset \prod_{\sigma \in [0,1]} C_0(\weights \times T, \mathbb{K}(\H))^W,  
$$ 
again taking into account the canonical rescaling $ \mathfrak{t}_{q^\sigma} \cong \mathfrak{t}/X_*(T) \cong T $ 
for all $ \sigma > 0 $. 

\begin{prop} \label{qassemblyfieldcontinuity}
The image of $ C^*_\red({\bf G}) $ under the above inclusion map is contained in $ C_0([0,1] \times \weights \times T, \mathbb{K}(\H))^W $. 
\end{prop} 

\proof We need to show that the sections $ \omega^{\alpha}_{ij} \bowtie u^\beta_{kl} $ 
in $ C^*_\red(G_{q^\sigma}) \subset C_0(\weights \times T, \mathbb{K}(\H))^W $ depend continuously on $ \sigma $. 
Since the images of these elements in the natural representations on $ \H_{\mu, \lambda} $ are finite rank operators, 
it is enough to show that they are strongly continuous in $ \sigma $. For the operators $ \omega^\alpha_{ij} $ this is obvious because their 
action is constant across the interval. For the operators $ u^\beta_{kl} $ strong continuity follows by inspecting the explicit formulas for 
the natural representations and the fact that multiplication of matrix elements in $ \CF^\infty(K_{q^\sigma}) $ depends continuously on $ \sigma \in [0,1] $. \qed 

\begin{theorem} \label{constructquantumassembly}
The algebra $ C^*_\red({\bf G}) $ is a continuous field of $ C^* $-algebras over $ [0,1] $ with fibres $ C^*_\red({\bf G})_\sigma = C^*_\red(G_{q^\sigma}) $, 
trivial away from zero. 
\end{theorem} 

\proof It follows immediately from Proposition \ref{qassemblyfieldcontinuity} that 
$ C^*_\red({\bf G}) $ is a continuous field of $ C^* $-algebras over $ [0,1] $, and that the restriction of this field to the open interval $ (0,1] $ 
agrees with the trivial field $ C_0((0,1] \times \weights \times T, \mathbb{K}(\H))^W $. This shows in particular that the fibre
$ C^*_\red({\bf G})_\sigma $ of $ C^*_\red({\bf G}) $ for $ \sigma > 0 $ is equal to $ C^*_\red(G_{q^\sigma}) $. 

To compute the fibre at $ \sigma = 0 $ remark that $ C^*_\red(G_1) \cong K \ltimes C(K) $ is a full crossed product. Therefore the 
natural quotient homomorphism $ C^*_\red({\bf G})_0 \rightarrow C^*_\red(G_{q^0}) = C^*_\red(G_1) $ splits. Moreover, the image of the splitting homomorphism 
$ C^*_\red(G_1) \rightarrow C^*_\red({\bf G})_0 $ is dense. In other words, the fibre of $ C^*_\red({\bf G}) $ 
at $ \sigma = 0 $ identifies canonically with $ C^*_\red(G_1) $. \qed 

The continuous field obtained in Theorem \ref{constructquantumassembly} will be called the \emph{quantum assembly field} for the group $ G $. 

Let us also define a continuous field of $ C^* $-algebras $ A^Q $ over $ [0,1] $ by setting 
$$
A^Q = \{F \in C_0([0,1] \times \weights \times T, \mathbb{K}(\H))^W \mid F_0(t) \in \KH(L^2(K))^{K_t} \text{ for all } t \in T \}. 
$$
Here $ F_0 \in C_0(\weights \times T, \mathbb{K}(\H))^W $ is the evaluation of $ F $ at $ 0 $ in the first variable. 
This field is clearly trivial away from $ 0 $, and has fibres $ A^Q_\sigma \cong C^*_\red(G_{q^\sigma}) $ due to Theorem \ref{qcmgplancherelform}
and Theorem \ref{reducedcstarstructurequantum}. From our above considerations we immediately obtain the following result. 

\begin{prop} \label{qbcfieldiso}
The quantum assembly field is isomorphic to the continuous field $ A^Q $ defined above. 
\end{prop} 

\proof Theorem \ref{constructquantumassembly} shows that $ C^*_\red({\bf G}) $ embeds into $ A^Q $, and it is straightforward to 
check directly that this embedding is surjective. \qed

\subsection{The quantum assembly map} 

In this subsection we show that the quantum assembly field induces an isomorphism 
$$ 
\mu_q: K_*(K \ltimes C(K)) \rightarrow K_*(C^*_\red(G_q)),  
$$
following the approach in \cite{Higsonmackeyanalogyktheory} for the classical Baum-Connes assembly map. 

Recall the following Lemma from \cite{Higsonmackeyanalogyktheory}. For the convenience of the reader we include a proof. 

\begin{lemma} \label{Higsonlemma}
Let $ B $ be a $ C^* $-algebra and let $ p \in M(B) $ be a projection. If $ p $ acts as a rank one projection in all irreducible representations 
of $ B $, then $ p B p $ is a commutative $ C^* $-algebra which is Morita equivalent to $ B $. 
\end{lemma} 

\proof Consider the representation $ \pi: B \rightarrow \bigoplus_{\omega \in PS(B)} \LH(\H_\omega) $, where $ \omega $ 
ranges over the set $ PS(B) $ of pure states of $ B $ and $ \H_\omega $ is the GNS-representation associated with $ \omega $. 
Then $ \pi $ is injective, and by assumption the image of $ p $ in $ \LH(\H_\omega) $ is a rank one projection for all $ \omega \in PS(B) $. 
It follows that $ pBp \cong \pi(p B p) = \pi(p) \pi(B) \pi(p) $ is a commutative $ C^* $-algebra. 

We claim that $ B = B p B $. Indeed, if $ B p B \subset B $ is a proper ideal then there exists an irreducible representation $ \rho: B \rightarrow \LH(\H) $ 
which vanishes on $ BpB $. This implies $ \rho(p) = 0 $ which contradicts our assumption that $ \rho(p) $ must be a rank one projection. 

It follows that the Hilbert $ B $-module $ \E = p B $ with the structure inherited from $ B $ viewed as a Hilbert $ B $-module over itself
implements a Morita equivalence between $ B $ and $ p B p $, see for instance example 3.6 in \cite{RaeburnWilliamsmorita}. \qed 

In particular, under the assumptions of Lemma \ref{Higsonlemma} the spectrum of $ B $ agrees with the spectrum of the commutative $ C^* $-algebra $ p B p $. 

Due to Theorem \ref{constructquantumassembly} the quantum assembly field $ C^*_\red({\bf G}) $ is trivial away from zero. Hence it induces 
a map in $ K $-theory 
$$
K_*(K \ltimes C(K)) = K_*(C^*_\red({\bf G})_0) \rightarrow K_*(C^*_\red({\bf G})_1) = K_*(C^*_\red(G_q))
$$
as in the classical case. We will call this map the \emph{quantum assembly map} for $ G $. 

\begin{theorem} \label{BCquantum}
Let $ G_q $ be a semisimple complex quantum group. The quantum assembly map 
$$
\mu_q: K_*(K \ltimes C(K)) \rightarrow K_*(C^*_\red(G_q))
$$
is an isomorphism. 
\end{theorem} 
\proof We follow closely the argument given by Higson in \cite{Higsonmackeyanalogyktheory} for the classical assembly map. 

Fix a labelling $ \weights^+ = \{\mu_1, \mu_2, \mu_3, \dots \} $ of the set of dominant integral weights such that 
if $ \mu_i \leq \mu_j $ with respect to the natural order on weights then $ i \leq j $. 

For each $ n \in \mathbb{N} $ we consider the projections $ p_n \in C^*_\red({\bf G}) = Q $ onto the highest weight vector of the 
minimal $ K $-type $ \mu_n $ and set 
$$
J_n = Q p_n Q.  
$$
In this way we obtain a filtration $ (\F^n(Q))_{n \in \mathbb{N}} $ of $ Q $ by ideals where 
$$
\F^n(Q) = \sum_{1 \leq m \leq n} J_n. 
$$
Let $ Q_1 = \F^1(Q) $ and 
$$
Q_n = \F^n(Q)/\F^{n - 1}(Q) 
$$
for $ n > 1 $ be the corresponding subquotients. 
The projections $ p_n \in Q_n $ satisfy the assumptions of Lemma \ref{Higsonlemma}, and according to Theorem \ref{reducedcstarstructurequantum} 
and Theorem \ref{qcmgplancherelform}, combined with Proposition \ref{qassemblyfieldcontinuity}, we get 
$$
p_n Q_n p_n \cong C([0,1], C(T)^{W_{\mu_n}}), 
$$
where $ W_{\mu_n} \subset W $ is the subgroup fixing $ \mu_n \in \weights^+ $. 

In the same way we obtain a filtration of $ D = C^*_\red(G_q) $ by considering the ideals 
$$
\F^n(D) = \sum_{1 \leq m \leq n} I_n
$$
where $ I_n = D p_n D $, and we will write $ D_n $ for the resulting subquotients. 
From the description of $ p_n Q_n p_n $ above and $ p_n D_n p_n \cong C(T)^{W_{\mu_n}} $ it follows that evaluation at $ 1 $ yields an 
isomorphism $ K_*(p_n Q_n p_n) \rightarrow K_*(p_n D_n p_n) $ in $ K $-theory. 

According to Lemma \ref{Higsonlemma}, in the commutative diagram 
$$
\xymatrix{
K_*(Q_n) \ar@{->}[r] \ar@{<-}[d] & K_*(D_n) \ar@{<-}[d] \\
K_*(p_n Q_n p_n) \ar@{->}[r] & K_*(p_n D_n p_n)
     }
$$
the vertical maps, induced by the canonical inclusions, are isomorphisms. Hence evaluation at $ 1 $ 
yields isomorphisms $ K_*(Q_n) \rightarrow K_*(D_n) $ for all $ n \in \mathbb{N} $. 
As a consequence, an iterated application of the $ 6 $-term exact sequence in $ K $-theory shows that 
evaluation at $ 1 $ induces isomorphisms $ K_*(\F^n(Q)) \rightarrow K_*(\F^n(D)) $ for all $ n \in \mathbb{N} $. 

Let us write $ \F^\infty(Q) $ and $ \F^\infty(D) $ for the union of the subalgebras $ \F^n(Q) $ and $ \F^n(D) $, respectively. Then 
$ \F^\infty(Q) \subset Q $ and $ \F^\infty(D) \subset D $ are dense, and 
continuity of $ K $-theory implies that evaluation at $ \sigma = 1 $ induces an isomorphism
$$
K_*(C^*_\red({\bf G})) = \varinjlim_{n \in \mathbb{N}} K_*(\F^n(Q)) \cong \varinjlim_{n \in \mathbb{N}} K_*(\F^n(D)) = K^*(C^*_\red(G_q)).  
$$

The map $ K_*(Q) = K_*(C^*_\red({\bf G})) \rightarrow K^*(C^*_\red(G_1)) = K_*(Q_0) $ induced 
by evaluation at $ \sigma = 0 $ is an isomorphism by the triviality of $ C^*_\red({\bf G}) $ away from $ 0 $. Hence the commutative diagram 
\begin{align*}
\xymatrix{
& K_*(C^*_\red({\bf G})) \ar[dl]_\cong \ar[dr]^\cong \\
K_*(C^*_\red(G_1)) \ar[rr]_{\mu_q} & & K_*(C^*_\red(G_q))
}
\end{align*}
yields the claim. \qed

\section{The deformation square} \label{secdefsquare}

We shall now assemble the constructions in previous sections and construct the \emph{deformation square} of a complex semisimple Lie group. 
This continuous field of $ C^* $-algebras encapsulates all the deformations we have studied so far. 

\subsection{Deformation squares and their $ K $-theory} 

In this subsection we introduce the abstract concept of a deformation square and discuss their basic $ K $-theoretic properties. 

\begin{definition} 
Let $ X = [0,1] \times [0,1] $ be the unit square and let $ A $ be a continuous field of $ C^* $-algebras over $ X $. 
We shall say that $ A $ is a deformation square if the following conditions hold. 
\begin{bnum}
\item[a)] The restriction of $ A $ to $ (0, 1] \times (0, 1] $ is a trivial field. 
\item[b)] The restriction of $ A $ to $ \{0\} \times [0,1] $ is trivial away from $ (0, 0) $. 
\item[c)] The restriction of $ A $ to $ [0,1] \times \{0\} $ is trivial away from $ (0, 0) $. 
\end{bnum} 
\end{definition} 

If $ A $ is a deformation square then the restriction of $ A $ to any boundary edge of $ [0,1] \times [0,1] $ yields a continuous field
of $ C^* $-algebras which is trivial away from a single point. The following result confirms that the corresponding induced maps in $ K $-theory 
behave as expected. 

\begin{prop} \label{Ktheorysquare}
Let $ A $ be a deformation square. Then one obtains a commutative diagram 
$$
\xymatrix{
K_*(A_{0,1}) \ar@{->}[r] \ar@{<-}[d] & K_*(A_{1,1}) \ar@{<-}[d] \\
K_*(A_{0,0}) \ar@{->}[r] & K_*(A_{1,0})
     }
$$
induced from the restriction of $ A $ to the boundary edges of $ X = [0,1] \times [0,1] $. 
\end{prop} 

\proof Using split exactness of $ K $-theory and Proposition 3.2 in \cite{Blancharddef} it suffices to consider the case that $ A $ is 
a unital deformation square. 

Let $ p_{0, 0} \in M_n(A_{0,0}) $ be a projection. In the same way as in our discussion in section \ref{secprelim} 
we can lift this to a positive element $ q $ in $ M_n(A) $. 
Since $ p_{0,0}^2 = p_{0, 0} $ we know that $ q(\sigma, \tau)^2 - q(\sigma, \tau) $ has small norm for $ (\sigma, \tau) $ near $ (0,0) $. 
Hence $ \frac{1}{2} $ is not in the spectrum of $ q(\sigma, \tau) $, 
and using functional calculus we obtain an element $ p \in A $ lifting $ p_{0,0} $ such that $ p(\sigma, \tau) $ is a projection 
provided $ 0 \leq \sigma, \tau \leq \epsilon $ for some $ \epsilon > 0 $. Using triviality of the field along the boundary and the interior, we can 
actually assume that $ p(\sigma, \tau) $ is a projection for all $ (\sigma, \tau) \in X $. By the very construction of the induced maps in $ K $-theory, 
both compositions of the above diagram map $ [p_{0,0}] $ to $ [p(1,1)] \in K_0(A_{1,1}) $. We conclude that the diagram for $ K_0 $ 
is commutative. 

In a similar way one proceeds for $ K_1 $. If $ u_{0, 0} \in M_n(A_{0,0}) $ is a unitary and $ u \in M_n(A) $ a lift of $ u_{0,0} $,  
then since $ u(\sigma, \tau) u(\sigma, \tau)^* - 1 $ and $ u(\sigma, \tau)^* u(\sigma, \tau) - 1 $ have small norm near $ 0 $, the 
elements $ u(\sigma, \tau) $ will be invertible provided $ 0 \leq \sigma, \tau \leq \epsilon $ for some $ \epsilon > 0 $. Using triviality of the 
field along the boundary and the interior we can in fact assume that $ u(\sigma, \tau) $ is invertible for all $ (\sigma, \tau) \in X $. Again, by 
construction of the induced maps in $ K $-theory, both compositions in the above diagram map $ [u_{0,0}] \in K_1(A_{0,0}) $ 
to $ [u(1,1)] \in K_1(A_{1,1}) $. \qed

\subsection{The classical assembly field revisited}

In this subsection we shall give an alternative description of the Baum-Connes assembly field for the classical group $ G $. 
More precisely, we shall prove that this field agrees with the continuous field 
$$
A^C = \{F \in C_0([0,1] \times \weights \times \mathfrak{t}, \mathbb{K}(\H))^W \mid F_0 \in A^L_0 \}
$$
over $ [0,1] $, where $ A^L_0 $ is the $ C^* $-algebra in Theorem \ref{cmgplancherelform} and $ F_0 \in C_0(\weights \times \mathfrak{t}, \mathbb{K}(\H))^W $ 
is the evaluation of $ F $ at $ 0 \in [0,1] $.
The field $ A^C $ is clearly trivial away from $ 0 $, and has fibres $ A^C_0 \cong K \ltimes C_0(\mathfrak{k}) $ 
and $ A^C_\sigma \cong C^*_\red(G) $ for $ \sigma > 0 $ according to Theorem \ref{cmgplancherelform} and Theorem \ref{reducedcstarstructure}, respectively. 

Recall that the group $ G $ admits the Iwasawa decomposition $ G = KAN $. That is, any element $ g \in G $ can be uniquely decomposed as 
$ g = k b $ for $ k \in K $ and $ b \in AN $. In particular, for any $ b \in AN $ and $ k \in K $ there exist unique elements $ b \hit k \in K $ 
and $ b \hitby k \in AN $ 
such that 
$$ 
bk = (b \hit k)(b \hitby k).  
$$
This leads to a left action of $ AN $ on $ K $ and a right action of $ K $ on $ AN $, respectively, which are called the left and right dressing actions. 
For more information see \cite{LuWeinsteinpoisson} and section 1.5 in \cite{CPbook}. 

\begin{theorem} \label{bcfieldiso}
The continuous field $ A^C $ defined above is isomorphic to the Baum-Connes assembly field for the group $ G $. 
\end{theorem}

\proof We shall verify that a set of generating sections of the assembly field define elements of the $ C^* $-algebra $ A^C $. For this we need to show that 
these sections are norm-continuous once we identify all fibres with the corresponding subalgebras of $ C_0(\weights \times \mathfrak{t}, \mathbb{K}(\H))^W $. 

Let us first recall that the strong-* topology on bounded subsets of $ \LH(\H) $ coincides with the strict topology. 
In particular, for a strong-* convergent net $ (T_i)_{i \in I} $ of uniformly bounded operators and a compact operator $ S $, the net $ (S T_i)_{i \in I} $ 
is norm convergent. 

Now let $ (f_\sigma)_{\sigma \in [0,1]} $ be a generating continuous section of the assembly field obtained from a 
function $ f \in C^\infty_c(K \times \mathfrak{an}) \subset C^*(K \ltimes \mathfrak{k}^*) $ 
which is invariant under convolution with some finite rank projection in $ \DF(K) $ from the left and right. 
According to our above remarks, it suffices to show that the operators $ \pi_{\mu, X}(f_\sigma) \in \KH(L^2(\E_\mu)) $ depend strong-* continuously 
on $ \sigma \in [0,1] $ and $ X \in \mathfrak{t} $ for all $ \mu \in \weights $. 

Let us compute $ \pi_{\mu, X}(f_\sigma)(\xi) $ for a generating section $ (f_\sigma)_{\sigma \in [0,1]} $ associated with $ f \in C^\infty_c(K \times \mathfrak{an}) $ 
as above and $ \xi \in L^2(\E_\mu) $. That is, we consider 
$$
f_0(k, a + n) = f(k, a + n)
$$
and 
$$
f_\sigma(k, \exp(a) \exp(n)) = f(k, \sigma^{-1} a + \sigma^{-1} n), \quad \sigma > 0.  
$$
Let $ \S(\mathfrak{k}) $ be the Schwartz space of $ \mathfrak{k} $ and $ \F: C^\infty_c(\mathfrak{an}) \rightarrow \S(\mathfrak{k}) $ the Fourier transform 
given by 
$$
\F(h)(Y) = \int_{\mathfrak{an}} e^{-i \Im(Y, a + n)} h(a + n) da dn. 
$$
Here $ \Im(\;, \;) $ is the imaginary part of the invariant bilinear form $ ( \;,\;) $ on $ \mathfrak{g} $. We note that 
this form induces a linear isomorphism between $ \mathfrak{an} $ and the real dual space $ \mathfrak{k}^* = \Hom(\mathfrak{k}, \mathbb{R}) $. 

Let us also write $ \F_{\mathfrak{an}}(f) $ for the Fourier transform of $ f $ in the $ \mathfrak{an} $-direction and denote by $ s \cdot X $ the 
adjoint action of $ s \in K $ on $ X \in \mathfrak{t} \subset \mathfrak{k} $. Then we have 
\begin{align*}
\pi_{\mu, X}(f_0)(\xi)(r) &= \int_{K} \F_{\mathfrak{an}}(f)(k, k^{-1} r \cdot X) \xi(k^{-1} r) dk \\
&= \int_{K \times \mathfrak{an}} f(k, a + n) e^{- (k^{-1}r \cdot X, a)} \xi(k^{-1} r) dk da dn. 
\end{align*}
The left Haar measure of $ G = KAN $ can be written as $ dx = e^{2(\rho, a)} dk da dn $ where $ dk $ is the Haar measure of $ K $ and 
$ da, dn $ is the Lebesgue measure on $ \mathfrak{a} \cong A $ and $ \mathfrak{n} \cong N $, respectively, see Proposition 8.43 in \cite{Knappbeyond}. 
Under the identification of the fibre algebra of the assembly field at $ \sigma > 0 $ with $ C^*_\red(G) $ we have to rescale the $ \mathfrak{an} $-part 
of this measure with $ \sigma^{-1} $. 
Using our notation $ b \hit k $ and $ b \hitby k $ for the dressing actions of $ b \in AN $ on $ k \in K $ 
and vice versa, we therefore obtain
\begin{align*}
&\pi_{\mu, X}(f_\sigma)(\xi)(r) \\
&= \int_{K \times \mathfrak{an}} f_\sigma(k, \exp(a) \exp(n)) \xi(\exp(-n) \exp(-a) k^{-1} r) e^{(2 \rho, a)} dk d(\sigma^{-1}a) d(\sigma^{-1}n) \\
&= \int_{K \times \mathfrak{an}} f(k, \sigma^{-1} a + \sigma^{-1} n) \xi(\exp(-n) \exp(-a) k^{-1} r) e^{(2 \rho, a)} dk d(\sigma^{-1}a) d(\sigma^{-1}n) \\
&= \int_{K \times \mathfrak{an}} f(k, \sigma^{-1} a + \sigma^{-1} n) \xi(\exp(-n) \exp(-a) \hit (k^{-1} r)) \\
&\qquad \times \chi_{\sigma^{-1} X + 2 \rho}(\exp(-n) \exp(-a) \hitby k^{-1} r) e^{(2 \rho, a)} dk d(\sigma^{-1}a) d(\sigma^{-1}n) \\
&= \int_{K \times \mathfrak{an}} f(k, a + n) \xi(\exp(-\sigma n) \exp(-\sigma a) \hit (k^{-1} r)) \\
&\qquad \times \chi_{\sigma^{-1} X + 2 \rho}(\exp(-\sigma a) \exp(-\sigma n) \hitby k^{-1} r) e^{(2 \rho, \sigma a)} dk da dn
\end{align*}
for $ \sigma > 0 $. Here $ \chi_{\sigma^{-1} X + 2 \rho} $ denotes the character of $ AN $ entering the definition of principal series representations as discussed 
in section \ref{seccg}. 

Let us compare the integrand of the last integral with the integrand in our formula for $ \pi_{\mu, X}(f_0)(\xi) $ above, 
splitting up the contribution coming from the character $ \chi_{\sigma^{-1} X + 2 \rho} $ as 
$$ 
\chi_{\sigma^{-1} X}(\exp(-\sigma a) \exp(-\sigma n) \hitby k^{-1} r) \chi_{2 \rho}(\exp(-\sigma a) \exp(-\sigma n) \hitby k^{-1} r). 
$$
Assume that $ f $ has support in $ K \times D $ where $ D \subset \mathfrak{an} $ is compact. The function 
$$ 
D \ni a + n \mapsto \chi_{2 \rho}(\exp(-\sigma a) \exp(-\sigma n) \hitby k^{-1} r) e^{(2 \rho, \sigma a)} 
$$ 
depends continuously on $ r, k \in K $ and converges uniformly to $ 1 $ as $ \sigma \rightarrow 0 $. 
Moreover, if $ C \subset \mathfrak{t} $ is compact we have 
\begin{align*}
\lim_{\sigma \rightarrow 0} &\chi_{\sigma^{-1} X}(\exp(-\sigma a) \exp(-\sigma n) \hitby k^{-1} r) \\
&= \chi_{X}\biggl(\frac{d}{d \sigma}(\exp(-\sigma a) \exp(-\sigma n) \hitby k^{-1} r)_{\sigma = 0}\biggr) \\
&= e^{-(X, a \cdot (k^{-1}r))} = e^{-(k^{-1}r \cdot X, a)},  
\end{align*}
uniformly in $ X \in C, r, k \in K $ and $ a + n \in D $. Here we use that the linearisation of the right dressing action of $ K $ at the 
identity of $ AN $ identifies with the right coadjoint action. 

Consider now $ \xi \in \Gamma(\E_\mu) \subset L^2(\E_\mu) \subset L^2(K) $. Then $ \xi $ is a continuous function on $ K $, 
and since the dressing action is continuous it follows that for $ \sigma \rightarrow 0 $ the function $ \xi(\exp(-\sigma n) \exp(-\sigma a) \hit (k^{-1} r)) $ 
converges uniformly to $ \xi(k^{-1} r) $ as a function of $ r, k \in K $ and $ a + n \in D $. 

In summary, we see that $ \pi_{\mu, X}(f_\sigma)(\xi) \rightarrow \pi_{\mu, X}(f_0)(\xi) $ in $ L^2 $-norm for all $ \xi \in \Gamma(\E_\mu) $, 
uniformly for $ X \in C $. Since $ \Gamma(\E_\mu) \subset L^2(\E_\mu) $ is dense we conclude 
that $ [0,1] \times \mathfrak{t} \ni (\sigma, X) \mapsto \pi_{\mu, X}(f_\sigma) \in \KH(L^2(\E_\mu)) $ is strongly continuous. 

In the same way one checks that $ [0,1] \times \mathfrak{t} \ni (\sigma, X) \mapsto \pi_{\mu, X}(f_\sigma)^* \in \KH(L^2(\E_\mu)) $ 
is strongly continuous; in fact, this follows from our above considerations because $ \pi_{\mu, X}(f_\sigma)^* = \pi_{\mu, X}((f^*)_\sigma) $ 
and $ f^* $ satisfies the same conditions as $ f $. 

Since the operators $ \pi_{\mu, X}(f_\sigma) $ are uniformly bounded we conclude that 
$ [0,1] \times \mathfrak{t} \ni (\sigma, X) \mapsto \pi_{\mu, X}(f_\sigma) \in \KH(L^2(\E_\mu)) $ is norm-continuous. 
We thus obtain a canonical $ C[0,1] $-linear $ * $-homomorphism from the assembly field into $ A^C $, and 
it is easy to check directly that this map is an isomorphism. \qed

\subsection{The deformation square of a complex semisimple Lie group} 

Let $ G $ be a complex semisimple Lie group. In this subsection we shall construct a deformation square over $ [0,1] \times [0,1] $ such that 
\begin{align*}
A_{0,0} &\cong K \ltimes C_0(\mathfrak{k}), \\
A_{\sigma, 0} &\cong C^*_\red(G), \quad \sigma > 0 \\ 
A_{0, \tau} &\cong K \ltimes C(K), \quad \tau > 0 \\ 
A_{\sigma, \tau} &\cong C^*_\red(G_{q^{\sigma \tau}}) \quad \text{ else. }
\end{align*} 
More precisely, let $ R \cong A^R $ be the quantization field for $ G $ and let $ L \cong A^L \subset A^R \cong R $ be the quantization field of the 
Cartan motion group. Setting
$$
A = \{F \in C([0,1], A^R) \mid F_0 \in A^L \}
$$
yields a continuous field of $ C^* $-algebras as desired. Here $ F_0 \in A^R $ denotes the evaluation of $ F $ at $ 0 $. 

\begin{theorem} \label{CGdeformationsquare}
The above construction yields a deformation square for any complex semisimple Lie group $ G $. 
\end{theorem} 

\proof Due to Theorem \ref{bcfieldiso} and Proposition \ref{qbcfieldiso} we see that the sections $ a = (a_{\sigma, \tau}) $ of this field satisfy 

\begin{bnum}
\item[a)] For each fixed $ \tau \in [0,1] $ the section $ \sigma \mapsto a_{\sigma, \tau} $ is a continuous section of the (quantum) assembly field. 
\item[b)] For each fixed $ \sigma \in [0,1] $ the section $ \tau \mapsto a_{\sigma, \tau} $ is a continuous section of the quantization field. 
\item[c)] On $ (0,1] \times (0,1] $ the section $ a $ can be identified with an element of the trivial continuous field with fibre $ C^*_\red(G_q) $. 
\end{bnum} 

This yields the claim. \qed 

It is natural to expect that the deformation square in Theorem \ref{CGdeformationsquare} is the total algebra of 
a locally compact quantum groupoid obtained from a continuous bundle of locally compact quantum groups, in the same way as the deformation picture 
of the classical Baum-Connes assembly map is associated with a bundle of locally compact groups. 
We shall not attempt to make this precise; for recent work relevant to this problem we refer to \cite{DCFquantumgln}.

\subsection{Classical and quantum Baum-Connes} 

We are now ready to assemble and summarize the results obtained above. 

The deformation picture of the classical Baum-Connes assembly map provides an isomorphism between the $ K $-theory of the 
Cartan motion group $ K \ltimes \mathfrak{k}^* $ and the $ K $-theory of the reduced group $ C^* $-algebra $ C^*_\red(G) $. 
According to Theorem \ref{BCquantum}, there is an analogous isomorphism between 
the $ K $-theory of the quantum Cartan motion group and the $ K $-theory of the reduced group $ C^* $-algebra $ C^*_\red(G_q) $. 

The $ K $-groups of the quantum Cartan motion group have been computed by Brylinski and Zhang \cite{BrylinskiZhangequivariantK}, and for the classical 
Cartan motion group the $ K $-groups are well-known \cite{BCH}, see also \cite{AdemGomezequivariantK} for a quite general treatment. 
Let us summarize the results as follows. 

\begin{theorem} \label{BZAG}
Let $ K $ be a simply connected compact Lie group of rank $ N $. Then we have 
$$
K_*(K \ltimes C_0(\mathfrak{k})) = 
\begin{cases} 
R(K) & * = \dim(K) \\ 
0 & \text{ else} 
\end{cases}
$$
where $ K $ acts on $ \mathfrak{k} $ via the adjoint action. 
Similarly, 
$$
K_*(K \ltimes C(K)) = \Lambda^* R(K)^N 
$$
for the crossed product with respect to the adjoint action of $ K $ on itself. 
\end{theorem} 

Here $ \Lambda^* R(K)^N $ denotes the exterior algebra of the free $ R(K) $-module $ R(K)^N $ over 
$ R(K) $. In particular, the rank of $ K_*(K \ltimes C(K)) $ as an $ R(K) $-module is $ 2^N $. 

Using the deformation square of $ G $ obtained in Theorem \ref{CGdeformationsquare} we see that the classical and quantum assembly maps 
are related as follows. 

\begin{theorem} \label{defsquaremain}
Let $ G $ be a simply connected complex semisimple Lie group with maximal compact subgroup $ K $ and let $ q \in (0,1) $. Then we obtain a commutative 
diagram 
$$
\xymatrix{
K_*(K \ltimes C(K)) \ar@{->}[r]^{\;\; \mu_q} \ar@{<-}[d] & K_*(C^*_\red(G_q)) \ar@{<-}[d] \\
K_*(K \ltimes C_0(\mathfrak{k})) \ar@{->}[r]^{\;\; \mu} & K_*(C^*_\red(G))
     }
$$
in $ K $-theory. Both horizontal maps are isomorphisms, and both vertical maps are split injective. 
\end{theorem}

\proof It follows from Theorem \ref{CGdeformationsquare} and Proposition \ref{Ktheorysquare} that we obtain a commutative
diagram as stated. The lower horizontal map is an isomorphism because the Baum-Connes conjecture holds for $ G $, 
and the upper horizontal map is an isomorphism by Theorem \ref{BCquantum}. 

It remains to verify the claim regarding the vertical maps. For this it suffices to consider the map 
$ K_*(K \ltimes C_0(\mathfrak{k})) \cong K^K_*(C_0(\mathfrak{k})) \rightarrow K^K_*(C(K)) \cong K_*(K \ltimes C(K)) $ 
on the left hand side. The generator of the free $ R(K) $-module $ K_*^K(C_0(\mathfrak{k})) \cong K^K_{* + n}(C_\tau(\mathfrak{k})) $ 
is given by the Bott element, here $ C_\tau(\mathfrak{k}) $ denotes the Clifford algebra bundle over $ \mathfrak{k} $ and $ n = \dim(\mathfrak{k}) $. 
The map $ K_*^K(C_0(\mathfrak{k})) \rightarrow K^K_*(C(K)) $ under consideration can be identified with 
$ \iota_* $, where $ \iota: C_0(\mathfrak{k}) \cong C_0(U) \rightarrow C(K) $ 
is the $ * $-homomorphism corresponding to the inclusion $ U \subset K $ of an open neighborhood $ U $ of $ e \in K $ which is $ K $-equivariantly 
diffeomorphic to $ \mathfrak{k} $. 

Let $ D_K \in KK^K_0(C_\tau(K), \mathbb{C}) $ be the Dirac element for $ K $, see \cite{Kasparov2}, \cite{Blackadarbook}. 
By slight abuse of notation we shall write again $ \iota $ for the $ * $-homomorphism $ C_\tau(\mathfrak{k}) \rightarrow C_\tau(K) $ corresponding to the 
inclusion $ \mathfrak{k} \cong U \subset K $. Then $ \iota^*(D_K) \in KK^K_0(C_\tau(\mathfrak{k}), \mathbb{C}) $ 
identifies with the Dirac element $ D_\mathfrak{k} \in KK^K_0(C_\tau(\mathfrak{k}), \mathbb{C}) $ for $ \mathfrak{k} $. 
It follows that $ R(K) \cong K^K_*(C_0(\mathfrak{k})) \rightarrow K^K_*(C(K)) $ is split injective, with  
splitting implemented by $ D_K $. \qed 

Let us conclude by taking a look at the special case $ G = SL(2, \mathbb{C}) $. In this case we have $ K = SU(2) $, and Theorem \ref{BZAG} reduces to 
\begin{align*}
K_*(K \ltimes C_0(\mathfrak{k})) &=  
\begin{cases} 
0 & * \text{ even} \\ 
R(K) & * \text{ odd} 
\end{cases} 
\\
K_*(K \ltimes C(K)) &=  
\begin{cases} 
R(K) & * \text{ even} \\ 
R(K) & * \text{ odd.} 
\end{cases} 
\end{align*}
Moreover, the natural map between these groups obtained from Theorem \ref{defsquaremain} corresponds to the identity in degree $ 1 $. 
Using the explicit descriptions given in Theorem \ref{reducedcstarstructure} and Theorem \ref{reducedcstarstructurequantum}, respectively, 
one can of course also calculate the $ K $-groups of $ C^*_\red(G) $ and $ C^*_\red(G_q) $ directly.

\bibliographystyle{hplain}

\bibliography{cvoigt}

\end{document}